\documentclass[10pt]{article}
\usepackage[english]{babel}
\usepackage[a4paper,margin=.8in]{geometry}
\usepackage[T1]{fontenc}
\usepackage[utf8]{inputenc}
\usepackage{lmodern}
\usepackage{microtype}
\usepackage{amsmath,amssymb,amsfonts,amsthm,mathtools}
\usepackage{graphicx}
\usepackage[margin=1cm]{caption}
\usepackage{subcaption}
\usepackage{color,fix-cm}
\usepackage{ dsfont }
\usepackage{verbatim}
\usepackage{array}
\usepackage{ wasysym }
\usepackage{tikz}
\usepackage{tkz-berge}
\usepackage{tkz-graph}
\usepackage{forloop}
\usepackage{ifthen}
\usepackage{tkz-euclide}
\usepackage{rotating}
\usepackage{enumitem}
\usepackage{makecell}
\usetikzlibrary{shapes.geometric}
\usetikzlibrary{decorations.shapes}
\usepackage{float}
\usepackage[mathlines,pagewise]{lineno}
\usepackage{aliascnt}
\usepackage{titling}
\usepackage[hang,flushmargin]{footmisc}
\usepackage[hidelinks]{hyperref}
\hypersetup{
  pdftitle={Proper edge-weightings of bipartite multigraphs},
  pdfauthor={Péter Madarasi and Márton Szeidemann},
  pdfsubject={Proper edge-weightings of bipartite multigraphs with two or three prescribed weights},
  pdfkeywords={proper edge-weighting, bipartite multigraph, odd multi-cactus, degree-constrained spanning tree, universal 3-edge-weightability}
}
\usepackage[nameinlink,noabbrev,capitalise,sort&compress]{cleveref}

\makeatletter
\edef\ReferenceTool@CleverefVersion{\csname ver@cleveref.sty\endcsname}%
\def\ReferenceTool@BrokenCleverefVersion{2018/03/27 v0.21.4 Intelligent cross-referencing}%
\ifx\ReferenceTool@CleverefVersion\ReferenceTool@BrokenCleverefVersion
  \let\cpageref\relax
  \DeclareRobustCommand{\cpageref}{%
    \@ifstar{\@crefstar{cpageref}}{\@cref{cpageref}}}%
  \let\Cpageref\relax
  \DeclareRobustCommand{\Cpageref}{%
    \@ifstar{\@crefstar{Cpageref}}{\@cref{Cpageref}}}%
  \providecommand*{\@setcpagerefrange}[3]{%
    \@@setcpagerefrange{#1}{#2}{cref}{#3}}%
  \providecommand*{\@setCpagerefrange}[3]{%
    \@@setcpagerefrange{#1}{#2}{Cref}{#3}}%
  \providecommand*{\@setlabelcpagerefrange}[3]{%
    \@@setcpagerefrange{#1}{#2}{labelcref}{#3}}%
\fi
\makeatother

\AtBeginDocument{%
}

\newcommand{\DeclareEquationCrefFormat}[1]{%
  \crefformat{#1}{##2\textup{(##1)}##3}%
  \Crefformat{#1}{##2\textup{(##1)}##3}%
  \crefrangeformat{#1}{##3\textup{(##1)}##4\nobreakdash--##5\textup{(##2)}##6}%
  \Crefrangeformat{#1}{##3\textup{(##1)}##4\nobreakdash--##5\textup{(##2)}##6}%
  \crefmultiformat{#1}%
    {##2\textup{(##1)}##3}%
    { and~##2\textup{(##1)}##3}%
    {, ##2\textup{(##1)}##3}%
    {, and~##2\textup{(##1)}##3}%
  \Crefmultiformat{#1}%
    {##2\textup{(##1)}##3}%
    { and~##2\textup{(##1)}##3}%
    {, ##2\textup{(##1)}##3}%
    {, and~##2\textup{(##1)}##3}%
  \crefrangemultiformat{#1}%
    {##3\textup{(##1)}##4\nobreakdash--##5\textup{(##2)}##6}%
    { and~##3\textup{(##1)}##4\nobreakdash--##5\textup{(##2)}##6}%
    {, ##3\textup{(##1)}##4\nobreakdash--##5\textup{(##2)}##6}%
    {, and~##3\textup{(##1)}##4\nobreakdash--##5\textup{(##2)}##6}%
  \Crefrangemultiformat{#1}%
    {##3\textup{(##1)}##4\nobreakdash--##5\textup{(##2)}##6}%
    { and~##3\textup{(##1)}##4\nobreakdash--##5\textup{(##2)}##6}%
    {, ##3\textup{(##1)}##4\nobreakdash--##5\textup{(##2)}##6}%
    {, and~##3\textup{(##1)}##4\nobreakdash--##5\textup{(##2)}##6}%
  \labelcrefformat{#1}{##2\textup{(##1)}##3}%
  \labelcrefrangeformat{#1}{##3\textup{(##1)}##4\nobreakdash--##5\textup{(##2)}##6}%
  \labelcrefmultiformat{#1}%
    {##2\textup{(##1)}##3}%
    { and~##2\textup{(##1)}##3}%
    {, ##2\textup{(##1)}##3}%
    {, and~##2\textup{(##1)}##3}%
  \labelcrefrangemultiformat{#1}%
    {##3\textup{(##1)}##4\nobreakdash--##5\textup{(##2)}##6}%
    { and~##3\textup{(##1)}##4\nobreakdash--##5\textup{(##2)}##6}%
    {, ##3\textup{(##1)}##4\nobreakdash--##5\textup{(##2)}##6}%
    {, and~##3\textup{(##1)}##4\nobreakdash--##5\textup{(##2)}##6}%
}
\DeclareEquationCrefFormat{equation}
\DeclareEquationCrefFormat{subequation}

\newcommand{\DeclareItemCrefFormat}[1]{%
  \crefformat{#1}{##2\textup{##1}##3}%
  \Crefformat{#1}{##2\textup{##1}##3}%
  \crefrangeformat{#1}{##3\textup{##1}##4\nobreakdash--##5\textup{##2}##6}%
  \Crefrangeformat{#1}{##3\textup{##1}##4\nobreakdash--##5\textup{##2}##6}%
  \crefmultiformat{#1}%
    {##2\textup{##1}##3}{ and~##2\textup{##1}##3}%
    {, ##2\textup{##1}##3}{, and~##2\textup{##1}##3}%
  \Crefmultiformat{#1}%
    {##2\textup{##1}##3}{ and~##2\textup{##1}##3}%
    {, ##2\textup{##1}##3}{, and~##2\textup{##1}##3}%
  \crefrangemultiformat{#1}%
    {##3\textup{##1}##4\nobreakdash--##5\textup{##2}##6}%
    { and~##3\textup{##1}##4\nobreakdash--##5\textup{##2}##6}%
    {, ##3\textup{##1}##4\nobreakdash--##5\textup{##2}##6}%
    {, and~##3\textup{##1}##4\nobreakdash--##5\textup{##2}##6}%
  \Crefrangemultiformat{#1}%
    {##3\textup{##1}##4\nobreakdash--##5\textup{##2}##6}%
    { and~##3\textup{##1}##4\nobreakdash--##5\textup{##2}##6}%
    {, ##3\textup{##1}##4\nobreakdash--##5\textup{##2}##6}%
    {, and~##3\textup{##1}##4\nobreakdash--##5\textup{##2}##6}%
  \labelcrefformat{#1}{##2\textup{##1}##3}%
  \labelcrefrangeformat{#1}{##3\textup{##1}##4\nobreakdash--##5\textup{##2}##6}%
  \labelcrefmultiformat{#1}%
    {##2\textup{##1}##3}{ and~##2\textup{##1}##3}%
    {, ##2\textup{##1}##3}{, and~##2\textup{##1}##3}%
  \labelcrefrangemultiformat{#1}%
    {##3\textup{##1}##4\nobreakdash--##5\textup{##2}##6}%
    { and~##3\textup{##1}##4\nobreakdash--##5\textup{##2}##6}%
    {, ##3\textup{##1}##4\nobreakdash--##5\textup{##2}##6}%
    {, and~##3\textup{##1}##4\nobreakdash--##5\textup{##2}##6}%
}
\DeclareItemCrefFormat{enumi}
\DeclareItemCrefFormat{enumii}
\DeclareItemCrefFormat{enumiii}
\DeclareItemCrefFormat{enumiv}

\theoremstyle{plain}
\newtheorem{theorem}{Theorem}

\newaliascnt{lemma}{theorem}
\newtheorem{lemma}[lemma]{Lemma}
\aliascntresetthe{lemma}

\newaliascnt{proposition}{theorem}
\newtheorem{proposition}[proposition]{Proposition}
\aliascntresetthe{proposition}

\newaliascnt{corollary}{theorem}
\newtheorem{corollary}[corollary]{Corollary}
\aliascntresetthe{corollary}

\newaliascnt{claim}{theorem}

\aliascntresetthe{claim}

\newaliascnt{conjecture}{theorem}

\aliascntresetthe{conjecture}

\theoremstyle{definition}
\newaliascnt{definition}{theorem}

\aliascntresetthe{definition}

\theoremstyle{remark}
\newaliascnt{remark}{theorem}
\newtheorem{remark}[remark]{Remark}
\aliascntresetthe{remark}

\newcommand{\DeclareNamedCrefType}[3]{%
  \crefname{#1}{#2}{#3}%
  \Crefname{#1}{#2}{#3}%
  \crefrangelabelformat{#1}{##3##1##4\nobreakdash--##5##2##6}%
}

\DeclareNamedCrefType{section}{Section}{Sections}
\DeclareNamedCrefType{subsection}{Section}{Sections}
\DeclareNamedCrefType{subsubsection}{Section}{Sections}
\DeclareNamedCrefType{appendix}{Appendix}{Appendices}
\DeclareNamedCrefType{figure}{Figure}{Figures}
\DeclareNamedCrefType{table}{Table}{Tables}
\DeclareNamedCrefType{footnote}{Footnote}{Footnotes}
\DeclareNamedCrefType{page}{Page}{Pages}

\DeclareNamedCrefType{theorem}{Theorem}{Theorems}
\DeclareNamedCrefType{lemma}{Lemma}{Lemmas}
\DeclareNamedCrefType{proposition}{Proposition}{Propositions}
\DeclareNamedCrefType{corollary}{Corollary}{Corollaries}
\DeclareNamedCrefType{claim}{Claim}{Claims}
\DeclareNamedCrefType{conjecture}{Conjecture}{Conjectures}
\DeclareNamedCrefType{definition}{Definition}{Definitions}
\DeclareNamedCrefType{remark}{Remark}{Remarks}

\makeatletter
\@onlypreamble\DeclareEquationCrefFormat
\@onlypreamble\DeclareItemCrefFormat
\@onlypreamble\DeclareNamedCrefType
\makeatother

\usepackage{xcolor}
\newcommand{\R}{\mathbb{R}}
\newcommand{\Z}{\mathbb{Z}}
\newcommand{\Q}{\mathbb{Q}}

\newcommand{\da}{d^{a}}
\newcommand{\db}{d^{b}}

\title{
  \vspace*{-25pt}
  Proper edge-weightings of bipartite multigraphs
}
\author{
  P\'eter Madarasi\thanks{HUN-REN Alfr\'ed R\'enyi Institute of Mathematics, Re\'altanoda u.\ 13--15, Budapest H-1053, Hungary; and Department of Operations Research, ELTE E\"otv\"os Lor\'and University, P\'azm\'any P. s.\ 1/c, Budapest H-1117, Hungary. E-mail: \texttt{madarasi@renyi.hu}, corresponding author}
  \and
  Márton Szeidemann\thanks{Department of Operations Research, ELTE E\"otv\"os Lor\'and University, P\'azm\'any P. s.\ 1/c, Budapest H-1117, Hungary. E-mail: \texttt{szeidi.marci.05@gmail.com}}
}
\date{\vspace{-.5cm}}

\begin{document}
\maketitle

\begin{abstract}
  Given a finite set $W$ of real numbers, a proper $W$-edge-weighting of a graph is an assignment of weights from $W$ to its edges such that adjacent vertices have distinct sums of incident edge weights.
  Thomassen, Wu, and Zhang showed that for positive integers $a<b$ with $a$ odd and $b$ even, the connected bipartite multigraphs without a proper $\{a,b\}$-edge-weighting are precisely the odd multi-cacti.
  We extend this characterization to all distinct real weights $a,b$ satisfying $ab(a+b)\neq0$.
  Using this characterization, we also prove that every connected bipartite multigraph with at least three vertices admits a proper $W$-edge-weighting for every three-element set $W\subseteq\R$, confirming the universal $3$-edge-weightability conjecture of Deng and Qiu for bipartite graphs.
  The proof of the two-weight characterization combines parity-based completion lemmas with a degree-constrained spanning-tree argument that makes a selected vertex a strict weighted local extremum among certain neighbors.
  For each fixed pair $a,b$ as above, it yields a polynomial-time algorithm constructing a proper $\{a,b\}$-edge-weighting or certifying that the graph is an odd multi-cactus.
  Both problems admit polynomial-time algorithms for rational weights, even when the edge multiplicities are encoded in binary.

\end{abstract}

\section{Introduction}\label{sec:introduction}

For a set $W$ of edge-weights, a $W$-edge-weighting of a graph $G$ is a map $\omega\colon E(G)\to W$.
We consider the induced weighted degrees
\[
  s_\omega(v)=\sum_{e\in E_G(v)}\omega(e),
\]
where $E_G(v)$ denotes the set of edges incident with $v$.
We call the edge-weighting $\omega$ \emph{proper} (or \emph{vertex-coloring}) if $s_\omega(u)\neq s_\omega(v)$ holds for every edge $uv \in E(G)$.
Equivalently, the induced weighted degrees form a proper vertex coloring of $G$.
We say that $G$ has the \emph{$W$-property} if it admits a proper $W$-edge-weighting.

We study proper edge-weightings of bipartite multigraphs with prescribed weight sets of size two or three.
Our first main result, \cref{thm:main}, characterizes the connected bipartite multigraphs with the $\{a,b\}$-property for distinct real weights $a,b$ satisfying $ab(a+b)\neq0$.
Our second main result, \cref{thm:universal}, establishes the existence of a proper edge-weighting from every three-element weight set for connected bipartite multigraphs of order at least three.

The notion of proper edge-weightings goes back to Karo\'nski, \L{}uczak, and Thomason, who formulated the 1-2-3 Conjecture, asking whether every simple graph without a component isomorphic to $K_2$ admits a proper $\{1,2,3\}$-edge-weighting~\cite{karonski2004edge}.
Keusch proved this conjecture~\cite{keusch2024solution123}.

A graph is called \emph{universally $3$-edge-weightable} if it has the $W$-property for every three-element set $W\subseteq\R$~\cite{deng2026regular}.
Deng and Qiu conjectured that every simple graph without a component isomorphic to $K_2$ is universally $3$-edge-weightable~\cite[Conjecture~7]{dengqiu2027uniformspan}.
Deng proved universal $3$-edge-weightability for regular simple graphs without a $K_2$ component~\cite[Corollary~2]{deng2026regular}.
The stronger list version of the 1-2-3 Conjecture, proposed by Bartnicki, Grytczuk, and Niwczyk~\cite{bartnicki2009weight}, asks whether every simple graph without a $K_2$ component admits a proper edge-weighting when each edge is given an arbitrary three-element list of real weights and its weight must be chosen from its own list.

Dudek and Wajc proved that deciding whether a given graph has the $\{1,2\}$-property is NP-complete, and they proved the same for $\{0,1\}$-edge-weightings~\cite{dudekwajc2011complexity}.
Thus, the two-weight problem is not only a special case of the general 1-2-3 problem; it has its own structural and algorithmic features.
Later, Madarasi and Simon extended their result: for any fixed distinct integers $a,b$, it is NP-complete to decide whether a simple cubic planar graph has the $\{a,b\}$-property~\cite{madarasi2026complexity}.

Bipartite graphs have been a central test class for two-weight results.
Chang, Lu, Wu, and Yu studied the minimum $k$ for which a graph has a vertex-coloring $\{1,\dots,k\}$-edge-weighting and determined this parameter for several classes, including trees and $r$-regular bipartite graphs with $r\geq3$~\cite{chang2011graphs}.
Lu, Yu, and Zhang proved several sufficient conditions for the existence of a vertex-coloring $\{1,2\}$-edge-weighting; in particular, every 3-connected bipartite graph has such a weighting~\cite{luyuzhang2011twoedgeweighting}.
Lu later proved that every 2-connected and 3-edge-connected bipartite graph has a vertex-coloring $W$-edge-weighting for each $W\in\{\{0,1\},\{1,2\}\}$~\cite{lu2016bipartite}.
These results give broad positive classes, but they do not identify the exact obstruction family for connected bipartite graphs.

The forbidden family for the $\{1,2\}$ case was found by Thomassen, Wu, and Zhang~\cite{thomassenwuzhang2016flows}.
They give two equivalent descriptions of odd multi-cacti, one using odd combinations and the other using cycle pasting; we use the latter description below.
They show that every odd multi-cactus has no proper $\{1,2\}$-edge-weighting, and their Theorem~5.1 proves the converse: every connected bipartite graph without a proper $\{1,2\}$-edge-weighting is an odd multi-cactus~\cite{thomassenwuzhang2016flows}.
Their paper also records in its concluding remarks that the odd-multi-cactus characterization remains valid when $a,b$ are integers with $0<a<b$, $a$ odd, and $b$ even~\cite{thomassenwuzhang2016flows}.

The behavior of the $\{0,1\}$ case confirms that the precise choice of the two weights is relevant.
Lyngsie proved that connected bridgeless bipartite graphs without the $\{0,1\}$-property are exactly odd multi-cacti, and also characterized the trees without a proper $\{0,1\}$-edge-weighting~\cite[Theorems~3 and~4]{lyngsie2018nsd01}.
Consequently, once bridges are allowed in the $\{0,1\}$ setting, odd multi-cacti are not the only connected bipartite obstructions.
Bensmail, Mc Inerney, and Lyngsie obtained another exact result in a different arithmetic regime: for an odd integer $a$, a 2-connected bipartite graph has a proper $\{a,a+2\}$-edge-weighting if and only if it is not an odd multi-cactus~\cite{bensmail2022oddab}.
They also characterized the trees without the $\{-1,1\}$-property, obtaining a non-trivial obstruction family different from those for $\{1,2\}$ and $\{0,1\}$~\cite[Theorem~2]{bensmail2022oddab}.
More recently, Madarasi characterized the trees without a proper $\{a,b\}$-edge-weighting for every pair of distinct real weights; in particular, when $ab(a+b)\neq0$, the only tree obstruction is $K_2$~\cite{madarasi2026trees}.
Since scaling both weights by a non-zero constant preserves properness, the tree obstructions other than $K_2$ for $\{0,1\}$ and $\{-1,1\}$ show that neither of the restrictions $ab\neq0$ and $a+b\neq0$ in~\cref{thm:main} can be omitted.

\subsection{Notation and terminology}

In the two-weight setting, $a$ and $b$ denote distinct real weights.
The convenient additional assumptions on the weights depend on the specific statements; these are locally stated in each section.
The cases where $a$ and $b$ have the same sign or opposite signs are handled separately in some proofs; we refer to these scenarios as the \emph{same-sign case} and the \emph{opposite-sign case}, respectively.
Introduce the notation $\delta = |b-a|$ for the absolute difference of the weights.
For integers $p,q$, write $[p,q]=\{i\in\Z:p\leq i\leq q\}$.

All graphs are finite, loopless multigraphs.
Degrees and edge cuts count multiplicity.
For a vertex set $X$ and an edge set $F$ of a graph $G$, let $G-X$ denote the graph obtained from $G$ by deleting the vertices in $X$ together with all incident edges, and let $G-F$ denote the graph obtained from $G$ by deleting the edges in $F$.
Let $\Delta_G(X)$ denote the set of edges with exactly one endpoint in $X$, put $d_G(X)=|\Delta_G(X)|$, and write $\sigma_G(X)$ for the number of connected components of $G-X$.
For a vertex $v$, let $E_G(v)$ denote the set of edges incident with $v$, and abbreviate $d_G(\{v\})$ to $d_G(v)$; thus $E_G(v)=\Delta_G(\{v\})$ and $d_G(v)=|E_G(v)|$.
If $F\subseteq E(G)$, then $d_F(v)=|F\cap E_G(v)|$.
Introduce $N_G(v)$ for the set of neighbors of $v$.

An \emph{edge-bundle} is an unordered pair of adjacent vertices together with all edges between them.
We often denote an edge-bundle by its endpoint pair $uv$; when a concrete edge is meant, this is stated explicitly.
An edge-bundle may have multiplicity one, so the word ``bundle'' does not imply that parallel edges are present.
The multiplicity of an edge-bundle is the number of edges between its endpoints.
An edge-bundle is called simple if it has multiplicity one and multiple otherwise.
An edge is called simple if it is the unique edge in its bundle.

For an edge-weighting $\omega$ of $G$, recall that $s_\omega(v)=\sum_{e\in E_G(v)}\omega(e)$.
When the weighting is fixed, we write $s(v)$ for $s_\omega(v)$.
We call an edge or edge-bundle $uv$ \emph{proper} if $s_\omega(u)\neq s_\omega(v)$; otherwise, we call it \emph{conflicting}.
The weighting $\omega$ is \emph{proper} if each edge is proper.
Since all edges in a bundle have the same endpoints, it is enough to check one inequality for each edge-bundle.
For an $\{a,b\}$-edge-weighting $\omega$, a vertex $v \in V(G)$, and a weight $\lambda \in \{a,b\}$, the number of edges incident with $v$ with weight $\lambda$ is called the \emph{$\lambda$-degree} of $v$, and we denote it by $d^{\lambda}_\omega(v)$.
Note that $d_G(v) = \da_\omega(v) + \db_\omega(v)$.
Furthermore, we have
\begin{equation}\label{eq:basic-sum}
  s_\omega(v)=a \da_\omega (v) + b \db_\omega(v) = a d_G(v)+ (b-a) \db_\omega (v).
\end{equation}

A graph is \emph{two-connected} if it is connected and has no cutvertex.
We allow two-vertex graphs to be two-connected under this convention.
We regard two parallel edges as a cycle of length two.
A \emph{block} is a maximal two-connected subgraph, and an \emph{endblock} is a block containing at most one cutvertex of the whole graph.
The \emph{block-cutvertex tree} of a connected graph is a tree whose nodes correspond to the blocks and cutvertices of the graph.
A block node and a cutvertex node are adjacent if and only if the cutvertex is contained in the block.
If $G$ is two-connected, then $G$ itself is its endblock.
In a connected graph that is not two-connected, the leaves of the block-cutvertex tree correspond precisely to the endblocks; all edges from an endblock to the rest of the graph are incident with its unique cutvertex.

\subsection{Our results}

For $m\geq1$, let $K_2^m$ denote the two-vertex loopless multigraph with exactly $m$ parallel edges between its vertices.
Following~\cite{thomassenwuzhang2016flows}, we call a graph an \emph{odd multi-cactus} if it is a graph $K_2^m$ or can be obtained by the following construction.
Start with a simple cycle of length at least six and congruent to $2$ modulo~$4$, with its edges colored alternately red and green.
Successively attach further such cycles by identifying a green edge of the new cycle with a green edge of the current graph, making no other identifications.
Finally, replace each green edge by an edge-bundle of arbitrary positive multiplicity.
Note that odd multi-cacti are bipartite.

\begin{theorem}\label{thm:main}
Let $a,b\in\R$ satisfy $a\neq b$ and $ab(a+b)\neq0$.
A connected bipartite loopless multigraph has a proper $\{a,b\}$-edge-weighting if and only if it is not an odd multi-cactus.
\end{theorem}

\begin{theorem}\label{thm:universal}
Every connected bipartite loopless multigraph with at least three vertices is universally $3$-edge-weightable.
\end{theorem}

The proof of~\cref{thm:main} first establishes that every odd multi-cactus is an obstruction for all pairs of distinct weights.
For the two-vertex base graphs, this is immediate; for the others, the pasted-cycle representation gives a direct proof and records the alternating bundle structure used in the subsequent reductions.
For the other direction of~\cref{thm:main}, suppose that a counterexample exists and choose one with the fewest vertices and then the fewest edges.
Completion lemmas for connected graphs of odd order allow us to prescribe weights on some designated edges and then choose the remaining weights while enforcing the inequalities that define a proper weighting.
Using these completion lemmas, we identify several local configurations in which prescribing weights on certain edges and extending the resulting partial weighting yields a proper weighting of the whole graph.
These local constructions and the minimality assumption yield structural restrictions that allow us to suppress the suspended paths of length three in a suitable endblock and select a vertex in the compressed graph.
A degree-constrained spanning-tree argument then gives the edge set needed to apply the completion lemmas, except in a specific cubic configuration, which is handled separately by switching the weights on a cycle.
For fixed weights $a,b$ satisfying the hypotheses of~\cref{thm:main}, the proof yields a polynomial-time algorithm for constructing a proper $\{a,b\}$-edge-weighting of a connected bipartite multigraph or determining that no such weighting exists, as explained in~\cref{sec:algorithm}.
The running time remains polynomial when bundle multiplicities are given in binary and the output weighting specifies in binary the number of edges receiving weight $b$ in each bundle.
If the weights $a,b$ are rational numbers given in binary, the running time is polynomial in the total encoding length of the graph and the weights.

In the proof of~\cref{thm:main}, our goal is to treat the same-sign and the opposite-sign cases as uniformly as possible.
We achieve this by defining certain concepts and notations differently in the two cases, while ensuring that the same statements hold for both.

The main ingredients of the proof of~\cref{thm:main} are the completion lemmas and the local-extremum argument used to construct a proper edge-weighting of the minimal counterexample.
The proof of the $\{1,2\}$ case strongly relies on the positivity and distinct parity of the weights~\cite{thomassenwuzhang2016flows}.
The completion lemmas extend the parity arguments used in earlier work~\cite{luyuzhang2011twoedgeweighting,thomassenwuzhang2016flows} to arbitrary distinct real weights.
In the same-sign case, a suitable local minimum or maximum of the degrees can be made a strict weighted extremum by prescribing weights on suitable edges incident with certain neighbors.
In the opposite-sign case, we instead seek a strict weighted local maximum, and a neighbor of larger degree may require several such prescribed edges, with the number depending on the difference in degrees.
We use Frank's degree-constrained spanning tree theorem to choose these edges while preserving the connectivity needed for the completion lemmas.
Compared with Lemma~4.3 in~\cite{thomassenwuzhang2016flows}, this gives the degree-dependent edge-deletion bounds needed in the opposite-sign case.

For~\cref{thm:universal}, the two-weight characterization leaves only odd multi-cacti and the weight sets of the form $\{-c,0,c\}$ with $c>0$ to consider.
We prove that odd multi-cacti with at least three vertices admit proper edge-weightings even from arbitrary three-element edge lists, and handle $\{-1,0,1\}$ by a parity argument.

The paper is organized as follows.
\Cref{sec:preliminaries} establishes some preliminary results, including the completion lemmas and list-coloring tools used to weight the alternating bundle-classes of paths and cycles.
\Cref{sec:omc} develops the pasted-cycle representation of odd multi-cacti, proves the obstruction direction, and studies small modifications of odd multi-cacti.
\Cref{sec:redsteps} gives several local configurations in which we can define a weighting and extend it to a proper weighting of the whole graph.
\Cref{sec:proof-completion} completes the proof of~\cref{thm:main} and describes the resulting algorithm.
Finally,~\cref{sec:universal} proves~\cref{thm:universal}.

\section{Preliminaries}\label{sec:preliminaries}
\subsection{Elementary observations}\label{subsec:elem-obs}
Fix arbitrary real numbers $a \neq b$.
This section states some elementary, although essential, observations about $\{a,b\}$-edge-weightings.
It is obvious from the definition that properness is invariant under scaling the weights by a non-zero real number.
That is, for any $\alpha \in \R \setminus \{0\}$, an $\{a,b\}$-edge-weighting $\omega$ is a proper $\{a,b\}$-edge-weighting if and only if $\alpha \omega$ is a proper $\{\alpha a, \alpha b\}$-edge-weighting.
Therefore, we may assume that one of the following three cases holds:
\begin{enumerate}
    \item One of the weights is zero.
    By scaling and relabeling, we may assume $a = 0$ and $b = 1$.
    The present paper does not seek a characterization for this weight set.
    \item The weights are non-zero with a rational ratio.
    By scaling, we may assume that $a,b \in \Z \setminus \{0\}$ are relatively prime.
    \item The weights are non-zero with an irrational ratio.
    By scaling, we may assume $a \in \Z \setminus \{0\}, b \in \R \setminus \Q$.
\end{enumerate}

The following lemma handles the case where the weights have an irrational ratio.
This lemma was stated in a slightly different form in~\cite[Proposition~2.1(ii)]{khatirinejad2012twoedgeweights} using linear independence in $\R$ over $\Q$.
\begin{lemma}\label{lem:irrational-ratio}
  Let $G=(V,E)$ be a graph with the $\{\alpha, \beta\}$-property for some real numbers $\alpha \neq \beta$, and let $a,b$ be non-zero real numbers with an irrational ratio.
  Then $G$ has the $\{a,b\}$-property.
\end{lemma}

\begin{proof}
    We may assume $a \in \Z \setminus \{0\}, b \in \R \setminus \Q$.
    Let $\omega$ be a proper $\{\alpha,\beta\}$-edge-weighting of $G$.
    From $\omega$, construct an $\{a,b\}$-edge-weighting $\omega'$ by replacing $\alpha$ with $a$ or $\beta$ with $b$ on each edge.
    Suppose, for a contradiction, that $uv$ is conflicting with respect to $\omega'$.
    Then, rearranging the equality, we have $(\da_{\omega'}(v) - \da_{\omega'}(u)) a + (\db_{\omega'}(v) - \db_{\omega'}(u)) b = 0$.
    Since $a$ is rational and $b$ is irrational, the $b$-degrees are forced to be equal, and this also forces the $a$-degrees to be equal as $a$ is non-zero.
    However, this implies that $uv$ is conflicting with respect to the original weighting $\omega$, contrary to properness of $\omega$.
\end{proof}

By~\cref{lem:irrational-ratio}, it suffices to construct concrete $\{a,b\}$-edge-weightings in the rational-ratio case.
In that case, we can prove a weaker, yet crucial, necessary condition for an edge to be conflicting.

\begin{lemma}\label{lem:conflict-necessary}
  Let $G=(V,E)$ be a graph and let $a, b \in \Z$ be distinct non-zero, relatively prime integers, and write $\delta = |b-a|$.
  If an edge $uv$ is conflicting with respect to an $\{a,b\}$-edge-weighting $\omega$ of $G$, then $d_G(u) \equiv d_G(v) \pmod{\delta}$ holds.
\end{lemma}

\begin{proof}
  Assume that $uv$ is conflicting, that is, $s_\omega(u) = s_\omega(v)$.
  By~\cref{eq:basic-sum}, we have $a(d_G(u)-d_G(v))+(b-a)(\db_\omega(u)-\db_\omega(v))=0$.
  Since $\delta=|b-a|$, this implies $a d_G(u) \equiv a d_G(v) \pmod{\delta}$.
  The integers $a$ and $\delta$ are relatively prime, as $a$ and $b$ are relatively prime.
  Consequently, we also have $d_G(u) \equiv d_G(v) \pmod{\delta}$.
\end{proof}

\subsection{Completion lemmas}

We begin with the classical parity criterion for the existence of a $T$-join.
Related formulations appear in the edge-weighting literature, for instance in~\cite{luyuzhang2011twoedgeweighting,thomassenwuzhang2016flows}.
For the sake of completeness, we include a short proof.

\begin{lemma}%
\label{lem:parity-subgraph}
    Let $H$ be a connected multigraph and let $p\colon V(H)\to\Z_2$.
    There is a set $F\subseteq E(H)$ with $d_F(v)\equiv p(v)\pmod2$ at every vertex if and only if $\sum_v p(v)\equiv0\pmod2$.
\end{lemma}

\begin{proof}
Necessity is the handshaking lemma.
For sufficiency, root a spanning tree and process its non-root vertices from the leaves upward, including the parent edge exactly when needed to correct the current parity.
The even total then gives the prescribed parity at the root.
\end{proof}

We continue with the completion lemma that forms the basis for several proofs.
We study when a partial $\{a,b\}$-edge-weighting defined only on a subset of edges of a connected bipartite graph $G$ can be extended to a proper edge-weighting of $G$.
Let $H$ be the subgraph of $G$ spanned by the edges without an initial weight.
The function $\varphi(v)$ will record the contributions from edges already assigned weights, or, in general, we will first consider $\varphi$ as any real-valued function on $V(H)$.
Let $\omega$ be an $\{a,b\}$-edge-weighting on $H$.
Then the resulting sum at $v \in V(H)$ in the original graph $G$ is
\[
S_{\varphi,\omega}(v) = \varphi(v) + s_\omega(v) = \varphi(v) + a d_H(v) + (b-a) \db_\omega(v).
\]
\emph{Switching} an edge means replacing its chosen weight by the other admissible weight; a cycle switch switches every edge of a cycle.
We first state the completion lemma in a general form; then we derive a version convenient for the applications.

\begin{lemma}\label{lem:completion-general}
Let $H$ be a connected loopless multigraph of odd order, let $N$ be a bipartite multigraph on $V(H)$, and let $\varphi\colon V(H)\to\R$.
Let $a,b$ be distinct real weights.
There is an $\{a,b\}$-edge-weighting $\omega$ of $H$ such that for every edge $uv\in E(N)$ and every integer $k$,
\begin{equation}\label{eq:state-target}
  S_{\varphi,\omega}(u)-S_{\varphi,\omega}(v)\neq 2k\delta.
\end{equation}
In particular, all the comparisons specified by $N$ are proper.
\end{lemma}

\begin{proof}
Put $\theta(v)=(\varphi(v)+a d_H(v))/(b-a)$ and let $K$ be the spanning subgraph of $N$ consisting of the edges $uv$ with $\theta(u)-\theta(v)\in\Z$.
On each component $Q$ of $K$, choose $\alpha_Q$ so that $\theta(v)-\alpha_Q\in\Z$ for every $v\in V(Q)$.
Choose a proper two-coloring $\gamma$ of $K$ using colors $0,1$ and, for every component $Q$ of $K$ and every $v\in V(Q)$, prescribe
\[
  \varrho(v)\equiv\gamma(v)-(\theta(v)-\alpha_Q)\pmod2.
\]
Since $|V(H)|$ is odd, $K$ has a component with an odd number of vertices.
Interchanging the two colors on this component also gives a proper two-coloring and changes the parity of $\sum_{v \in V(H)}\varrho(v)$, so we can enforce an even total and realize $\varrho$ with a subgraph of $b$-edges in $H$ by~\cref{lem:parity-subgraph}; this defines an $\{a,b\}$-edge-weighting $\omega$.
Since $\db_\omega(v)\equiv\varrho(v)\pmod2$ holds for each $v\in V(H)$, the normalized values $S_{\varphi,\omega}/(b-a)=\theta+\db_\omega$ differ by an odd integer across each edge of $K$.
On every other edge of $N$, the normalized difference is non-integral.
Both possibilities give~\cref{eq:state-target}, since $2(b-a)\Z=2\delta\Z$.
\end{proof}

\begin{remark}
    Every edge-weighting $\omega$ whose $b$-degrees match the constructed parity prescription $\varrho$ in the proof of~\cref{lem:completion-general} satisfies the conclusion of the lemma.
\end{remark}

\begin{corollary}\label{cor:completion-applied}
    Let $G$ be a connected loopless bipartite multigraph, and let $H$ be a connected subgraph of $G$ of odd order.
    Let $a,b$ be distinct real weights.
    Let $\omega_0\colon E(G) \setminus E(H) \to \{a,b\}$ be a partial $\{a,b\}$-edge-weighting defined on the edges outside $H$.
    Then $\omega_0$ can be extended to $E(H)$ such that the difference of the total weighted degrees across each edge induced by $V(H)$ is not an integral multiple of $2 \delta$.
    In particular, each edge induced by $V(H)$ is proper, and switching any cycle in $G$ preserves this property.
\end{corollary}

\begin{proof}
    Apply~\cref{lem:completion-general} with $H$, $N = G[V(H)]$, and $\varphi$, where $\varphi(v) = \sum_{e\in E_G(v)\setminus E(H)}\omega_0(e)$.
    This gives an $\{a,b\}$-edge-weighting $\omega_1$ of $H$ such that $(\varphi(u) + s_{\omega_1}(u)) - (\varphi(v) + s_{\omega_1}(v))$ is not an integral multiple of $2 \delta$ for any edge $uv \in E(G)$ with $u,v \in V(H)$.
    Together, $\omega_0$ and $\omega_1$ form an $\{a,b\}$-edge-weighting $\omega$ of $G$.
    By construction, we have $s_\omega(w) = \varphi(w) + s_{\omega_1}(w)$ for any vertex $w \in V(H)$.
    Hence, the condition above means that $s_\omega(u) - s_\omega(v)$ is not an integral multiple of $2 \delta$.

    The properness of each edge induced by $V(H)$ is obvious.
    If we switch an edge, then the weighted degrees of the endpoints increase or decrease by $\delta$.
    Since the degrees in any cycle are even ($2$ for the vertices on the cycle, $0$ for the others), a cycle switch changes each weighted degree by an integral multiple of $2\delta$.
    Therefore, the difference of weighted degrees across any edge induced by $V(H)$ remains distinct from each integral multiple of $2\delta$.
\end{proof}

\begin{corollary}\label{cor:odd-order}
Every connected bipartite multigraph of odd order has a proper $\{a,b\}$-edge-weighting for any distinct real weights $a \neq b$.
\end{corollary}

\begin{proof}
Apply~\cref{cor:completion-applied} with $G=H$ and $\omega_0$ being the empty function.
\end{proof}

\subsection{List-coloring tools for paths and cycles}\label{sec:paths-cycles}

We first prove two elementary list-coloring lemmas for cycles and paths, and then apply the first one to edge-weightings of bipartite multigraphs.

\begin{lemma}\label{lem:list-cycle}
Let $n\geq2$, and let $L_1,\dots,L_n$ be finite lists with $|L_i|\geq2$.
Indices in this lemma are taken modulo~$n$.
There are values $x_i\in L_i$ such that $x_i\neq x_{i+1}$ for every $i\in[1,n]$, unless $n$ is odd and there is a two-element set $C$ such that $L_i=C$ for every $i\in[1,n]$.
\end{lemma}

\begin{proof}
If some list, say $L_1$, has at least three elements, then choose $x_2,\dots,x_n$ successively, with consecutive values different.
Then choose $x_1\in L_1$ different from both $x_2$ and $x_n$, which excludes at most two values.
If all lists have size two but are not equal, then relabel so that $L_1\neq L_n$, choose $x_1\in L_1\setminus L_n$, and continue greedily to $x_n$.
Then $x_n\neq x_1$ holds automatically.
When all lists are equal to the same two-element set, the two elements must be chosen alternately; therefore, a proper choice exists exactly when $n$ is even.
\end{proof}

\begin{lemma}\label{lem:list-path}
Let $n\geq1$, and let $L_1,\dots,L_n$ be lists such that $|L_1|\geq1$, $|L_n|\geq1$, and $|L_i|\geq2$ for $i\in[2,n-1]$.
We seek values $x_i\in L_i$ such that $x_i\neq x_{i+1}$ for $i\in[1,n-1]$.
Such a proper choice exists if $n=1$.
For $n\geq2$, such a proper choice exists in either of the following cases:
\begin{enumerate}[label=(\roman*),ref=(\roman*)]
\item\label{lem:list-path:large-list} $|L_1|\geq2$, or $|L_n|\geq2$, or $|L_i|\geq3$ for some $i\in[2,n-1]$;
\item\label{lem:list-path:alternating} $n$ is even, there is a two-element set $C$ such that $L_i=C$ for every $i\in[2,n-1]$, and there exist distinct values $c_L\in L_1\cap C$ and $c_R\in L_n\cap C$.
\end{enumerate}
\end{lemma}

\begin{proof}
The case $n=1$ is immediate, so assume $n\geq2$.

To prove~\cref{lem:list-path:large-list}, first suppose $|L_1|\geq2$.
Choose $x_n\in L_n$ arbitrarily.
For $i\in[2,n-1]$ in decreasing order, choose $x_i\in L_i\setminus\{x_{i+1}\}$.
Then choose $x_1\in L_1\setminus\{x_2\}$.
The case $|L_n|\geq2$ is symmetric.

It remains to prove this case when $|L_j|\geq3$ for some $j\in[2,n-1]$.
Choose $x_1,\dots,x_{j-1}$ greedily from left to right and choose $x_n,\dots,x_{j+1}$ greedily from right to left.
When choosing $x_j$, at most the two neighboring values are forbidden, so a value in $L_j$ remains.

To prove~\cref{lem:list-path:alternating}, choose distinct values $c_L\in L_1\cap C$ and $c_R\in L_n\cap C$.
Since $C$ has two elements and $n$ is even, the alternating sequence in $C$ that starts with $c_L$ ends with $c_R$.
This sequence is a proper choice.
\end{proof}

\begin{lemma}\label{lem:cycle-bundle}
Let $a,b$ be distinct real weights.
Let $G$ be a connected bipartite multigraph whose underlying simple graph is the simple cycle $C_{2r}$, where $r\geq2$, with edge-bundles $E_1,\dots,E_{2r}$ in cyclic order.
Then the graph $G$ has a proper $\{a,b\}$-edge-weighting unless $r$ is odd and one of the two alternating classes of bundles consists entirely of simple bundles.
Equivalently, the only obstruction in this cycle-bundle setting occurs for cycles of length $2$ modulo~$4$ with one alternating class of simple bundles.
\end{lemma}

\begin{proof}
Let $W_i$ be the total weight of $E_i$.
At the vertex between $E_{i-1}$ and $E_i$, the sum is $W_{i-1}+W_i$.
The weighted degrees at the two endpoints of $E_i$ are different exactly when $W_{i-1}\neq W_{i+1}$.
The odd-indexed and even-indexed bundles give two independent list-coloring problems on cycles of length $r$, where the list for $E_i$ is the set of possible total weights $\{|E_i| a+k(b-a):k\in[0,|E_i|]\}$.
Apply~\cref{lem:list-cycle}.
The exceptional case for an alternating class occurs precisely when $r$ is odd and every bundle in that class has multiplicity one.
\end{proof}

\section{Odd multi-cacti}\label{sec:omc}

\subsection{Pasted-cycle representation}

We now encode the structure of odd multi-cacti by a tree, in analogy to the block-cutvertex tree.
This representation gives a direct proof of the obstruction property and makes the local configurations used in the subsequent reductions explicit.

An \emph{alternately simple cycle} is a connected loopless multigraph whose underlying simple graph is a cycle of even length.
Its edge-bundles are colored alternately red and green, with every red bundle simple and every green bundle having arbitrary positive multiplicity.
Its length is the number of its edge-bundles.

A \emph{pasted-cycle representation} of a connected bipartite multigraph $G$ is specified by a tree $T$ with a fixed bipartition of its nodes into \emph{cycle nodes} and \emph{green-bundle nodes}, with at least one cycle node.
Assign a positive integer multiplicity label to each green-bundle node and an alternately simple cycle of length $2$ modulo~$4$ to each cycle node.
The neighbors of each cycle node are in bijection with the green bundles of its cycle, each having the multiplicity specified by the corresponding label.
To obtain $G$, start with disjoint copies of the assigned cycles and identify the bundles corresponding to each green-bundle node, matching endpoints so that the bipartitions agree and retaining their common multiplicity.
No other vertices belonging to distinct cycles are identified, and all red bundles remain simple.
In a pasted-cycle representation, each vertex belongs to the bundle of a unique green-bundle node, since the green bundles of each cycle form a matching covering its vertices and identifications occur only through those nodes.

\begin{figure}[h]
    \centering
    \begin{tikzpicture}[thick]
    \begin{scope}[shift={(3.5,0)}]
        \begin{scope}[shift = {(1.8,0.5)},rotate=192]
            \tikzset{VertexStyle/.style = {shape = circle, fill = black , minimum size = 7pt, inner sep=0pt}}
            \Vertex[x=0,y=0,L=$\phantom{}$, LabelOut=True]{U}
            \tikzset{VertexStyle/.style = {shape = circle, fill = green , minimum size = 5pt, inner sep=0pt}}
            \foreach \i in {1,...,2} {
                \pgfmathsetmacro{\angle}{(360/3)*\i}
                \pgfmathsetmacro{\x}{cos(\angle)}
                \pgfmathsetmacro{\y}{sin(\angle)}
                \Vertex[x=0.618*\x,y=0.618*\y,L=$\phantom{}$, LabelOut=True]{u\i}
                \draw[] (U) -- (u\i);
            }
            \node[below right] at (u1) {$2$};
            \node[shift = {(78:3mm)}] at (u2) {$1$};
        \end{scope}
        \begin{scope}[shift = {(-1.8,0.5)},rotate = -12]
            \tikzset{VertexStyle/.style = {shape = circle, fill = black , minimum size = 7pt, inner sep=0pt}}
            \Vertex[x=0,y=0,L=$\phantom{}$, LabelOut=True]{V}
            \tikzset{VertexStyle/.style = {shape = circle, fill = green , minimum size = 5pt, inner sep=0pt}}
            \foreach \i in {1,...,2} {
                \pgfmathsetmacro{\angle}{(360/3)*\i}
                \pgfmathsetmacro{\x}{cos(\angle)}
                \pgfmathsetmacro{\y}{sin(\angle)}
                \Vertex[x=0.618*\x,y=0.618*\y,L=$\phantom{}$, LabelOut=True]{v\i}
                \draw[] (V) -- (v\i);
            }
            \node[below left] at (v2) {$3$};
            \node[shift = {(102:3mm)}] at (v1) {$4$};
        \end{scope}
        \begin{scope}[shift = {(0,0)}, rotate = 18]
            \tikzset{VertexStyle/.style = {shape = circle, fill = black , minimum size = 7pt, inner sep=0pt}}
            \Vertex[x=0,y=0,L=$\phantom{}$, LabelOut=True]{W}
            \tikzset{VertexStyle/.style = {shape = circle, fill = green , minimum size = 5pt, inner sep=0pt}}
            \foreach \i in {1,...,5} {
                \pgfmathsetmacro{\angle}{(360/5)*\i}
                \pgfmathsetmacro{\x}{cos(\angle)}
                \pgfmathsetmacro{\y}{sin(\angle)}
                \Vertex[x=\x,y=\y,L=$\phantom{}$, LabelOut=True]{w\i}
                \draw[] (W) -- (w\i);
            }
        \end{scope}
        \draw[] (U) -- (w5);
        \draw[] (V) -- (w2);
        \node[above=1mm] at (w1) {$3$};
        \node[shift = {(78:3mm)}] at (w2) {$1$};
        \node[shift = {(102:3mm)}] at (w5) {$3$};
        \node[below left] at (w3) {$5$};
        \node[below right] at (w4) {$2$};
    \end{scope}

    \begin{scope}[shift={(-3.5,0)}]
        \begin{scope}[rotate = 360*3/10, shift = {(0,-1.487)}]
            \tikzset{VertexStyle/.style = {shape = circle, fill = black , minimum size = 3pt, inner sep=0pt}}
            \foreach \i in {1,...,6} {
                \pgfmathsetmacro{\angle}{(360/6)*\i}
                \pgfmathsetmacro{\x}{cos(\angle)}
                \pgfmathsetmacro{\y}{sin(\angle)}
                \Vertex[x=0.618*\x,y=0.618*\y,L=$\phantom{}$, LabelOut=True]{u\i}
            }
            \draw[green] (u1) -- (u2);
            \draw[green] (u1) to[bend left = 40] (u2);
            \draw[green] (u1) to[bend right = 40] (u2);
            \draw[red] (u2) -- (u3);
            \draw[green] (u3) to[bend left = 30] (u4);
            \draw[green] (u3) to[bend right = 30] (u4);
            \draw[red] (u4) -- (u5);
            \draw[green] (u5) -- (u6);
            \draw[red] (u6) -- (u1);
        \end{scope}
        \begin{scope}[rotate = -360*3/10, shift = {(0,-1.487)}]
            \tikzset{VertexStyle/.style = {shape = circle, fill = black , minimum size = 3pt, inner sep=0pt}}
            \foreach \i in {1,...,6} {
                \pgfmathsetmacro{\angle}{(360/6)*\i}
                \pgfmathsetmacro{\x}{cos(\angle)}
                \pgfmathsetmacro{\y}{sin(\angle)}
                \Vertex[x=0.618*\x,y=0.618*\y,L=$\phantom{}$, LabelOut=True]{v\i}
            }
            \draw[green] (v1) -- (v2);
            \draw[red] (v2) -- (v3);
            \draw[green] (v3) to[bend left = 45] (v4);
            \draw[green] (v3) to[bend left = 15] (v4);
            \draw[green] (v3) to[bend right = 15] (v4);
            \draw[green] (v3) to[bend right = 45] (v4);
            \draw[red] (v4) -- (v5);
            \draw[green] (v5) -- (v6);
            \draw[green] (v5) to[bend left = 40] (v6);
            \draw[green] (v5) to[bend right = 40] (v6);
            \draw[red] (v6) -- (v1);
        \end{scope}
        \begin{scope}[shift = {(0,0)}]
            \tikzset{VertexStyle/.style = {shape = circle, fill = black , minimum size = 3pt, inner sep=0pt}}
            \foreach \i in {1,...,10} {
                \pgfmathsetmacro{\angle}{(360/10)*\i}
                \pgfmathsetmacro{\x}{cos(\angle)}
                \pgfmathsetmacro{\y}{sin(\angle)}
                \Vertex[x=\x,y=\y,L=$\phantom{}$, LabelOut=True]{w\i}
            }
            \draw[red] (w1) -- (w2);
            \draw[green] (w2) -- (w3);
            \draw[green] (w2) to[bend left = 40] (w3);
            \draw[green] (w2) to[bend right = 40] (w3);
            \draw[red] (w3) -- (w4);
            \draw[green] (w4) -- (w5);
            \draw[red] (w5) -- (w6);
            \draw[green] (w6) -- (w7);
            \draw[green] (w6) to[bend left = 60] (w7);
            \draw[green] (w6) to[bend left = 30] (w7);
            \draw[green] (w6) to[bend right = 30] (w7);
            \draw[green] (w6) to[bend right = 60] (w7);
            \draw[red] (w7) -- (w8);
            \draw[green] (w8) to[bend left = 30] (w9);
            \draw[green] (w8) to[bend right = 30] (w9);
            \draw[red] (w9) -- (w10);
            \draw[green] (w10) -- (w1);
        \end{scope}
    \end{scope}

    \draw[<->] (-0.3,0) -- (0.3,0);
    \end{tikzpicture}
    \caption{Illustration of an odd multi-cactus and its incidence tree.}
    \label{fig:omc-fit}
\end{figure}
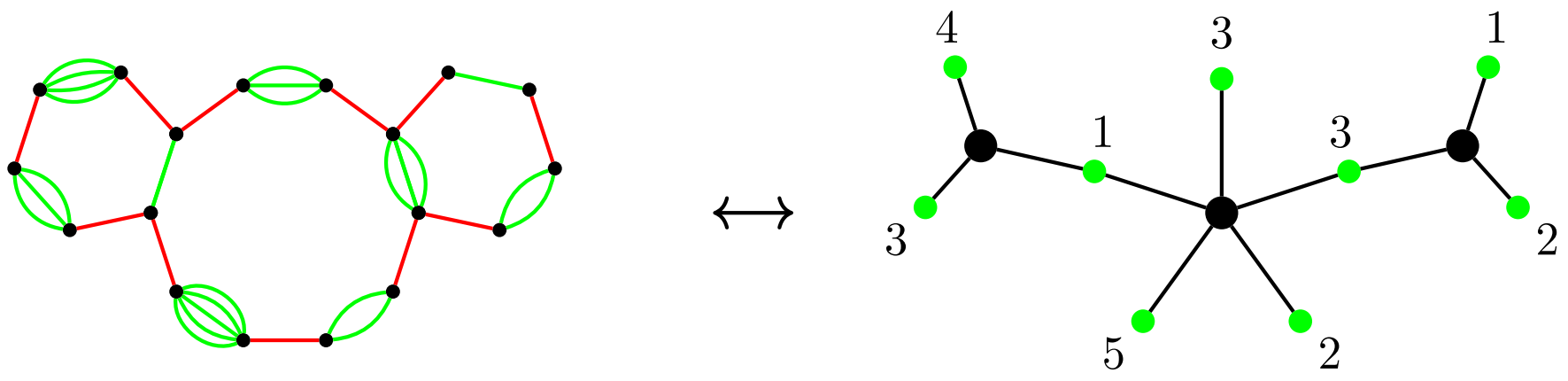

The image in $G$ of the alternately simple cycle assigned to a cycle node is the \emph{cycle subgraph} of that node.
We call $T$ the \emph{incidence tree}.
\Cref{fig:omc-fit} illustrates an incidence tree alongside the corresponding graph.
A green-bundle node of degree one in $T$ is \emph{private}; one of degree at least two is \emph{pasted}.
A green bundle is private if its green-bundle node has degree one and pasted if that node has degree at least two.
The exceptional graphs $K_2^m$ have no pasted-cycle representation.

\begin{proposition}
\label{prop:models}
The odd multi-cacti are exactly the graphs that admit a pasted-cycle representation, together with the exceptional graphs $K_2^m$.
\end{proposition}

\begin{proof}
The initial cycle gives an incidence tree consisting of a star with a central cycle node and one leaf for each green edge.
Each attachment adds a new star to the tree by identifying the green-bundle nodes corresponding to the common edge.
Assigning the final multiplicities to the green-bundle nodes gives a pasted-cycle representation.

Conversely, root the incidence tree at a cycle node and introduce the underlying simple cycles of the cycle subgraphs in nondecreasing order of the distances of their nodes from the root.
Attach each non-root cycle along the green edge corresponding to its parent green-bundle node, with no other identifications.
Assigning the multiplicities specified by the labels at the end yields an odd multi-cactus.
\end{proof}

We call a loopless multigraph $P$ an \emph{alternately simple path} if its underlying simple graph is a path of odd length and its edge-bundles are colored alternately red and green, starting with red, with each red bundle simple.
Its length is the number of its edge-bundles.
Deleting a green bundle from an alternately simple cycle yields an alternately simple path between its endpoints.

\begin{proposition}\label{prop:path-replacement}
Replacing an edge of a multiple bundle of an odd multi-cactus by an alternately simple path of length $4k+1$ for some positive integer $k$, with new internal vertices, yields an odd multi-cactus.
Conversely, every odd multi-cactus can be obtained from some $K_2^m$ by a possibly empty sequence of such replacements.
\end{proposition}

\begin{proof}
For $K_2^m$ with $m\geq2$, the result is an alternately simple cycle of length $4k+2$.
Otherwise, the multiple bundle is green in any pasted-cycle representation.
Decrease its label by one and attach to its node a new cycle node for the cycle formed by the remaining bundle and the path, with private green-bundle nodes labeled by the multiplicities of the green bundles of the path.

For the converse, the graphs $K_2^m$ require no replacements.
For any other odd multi-cactus, fix a pasted-cycle representation and root its incidence tree at a green-bundle node.
For each green-bundle node with multiplicity label $m$ and $t$ child cycle nodes, use multiplicity $m+t$ when its bundle first appears in the construction.
Start with the root bundle and introduce the cycle subgraphs in nondecreasing order of the distances of their nodes from the root.
For each cycle subgraph, replace an extra edge of its parent green bundle with the alternately simple path obtained by deleting that bundle from the cycle subgraph.
Use the augmented multiplicities for the other green bundles.
The internal vertices of the path are new by the tree structure.
Each replacement consumes one extra edge, so every green bundle remains non-empty and attains its prescribed multiplicity.
\end{proof}

\subsection{Obstructions and their small modifications}\label{subsec:omc-obs-mod}

Let $a,b$ be arbitrary distinct real weights.
In what follows, we establish that odd multi-cacti do not admit a proper $\{a,b\}$-edge-weighting; however, two kinds of small modifications of an odd multi-cactus result in a graph with the $\{a,b\}$-property.
The obstruction property for $\{1,2\}$ was established in~\cite{thomassenwuzhang2016flows}; we give a direct proof for arbitrary distinct real weights using the pasted-cycle representation.

\begin{lemma}\label{lem:omc-obstruction}
For any distinct real weights $a,b$, no odd multi-cactus has a proper $\{a,b\}$-edge-weighting.
\end{lemma}

\begin{proof}
For $K_2^m$, the two vertices are incident with the same edges and therefore have equal weighted degrees under every weighting.
Otherwise, by~\cref{prop:models}, fix a pasted-cycle representation of the odd multi-cactus $G$ with incidence tree $T$, and consider an arbitrary $\{a,b\}$-edge-weighting of $G$.
Orient each edge of $T$ toward its green-bundle node if the two red edges of the corresponding cycle subgraph incident with the endpoints of that green bundle have equal weights, and toward its cycle node otherwise.
Every cycle node has positive outdegree, since its odd number of red edges cannot alternate between the two weights in cyclic order.
Since every orientation of a tree is acyclic, some node has outdegree zero, and this node must be a green-bundle node.
The two endpoints of its bundle receive equal contributions from the bundle itself and from the two red edges of each incident cycle subgraph, so the bundle is conflicting.
Thus no $\{a,b\}$-edge-weighting of $G$ is proper.
\end{proof}

\begin{remark}\label{rem:mainthm-irratratio}
    By the known $\{1,2\}$-characterization and~\cref{lem:irrational-ratio,lem:omc-obstruction}, \cref{thm:main} holds whenever $a,b$ have an irrational ratio.
\end{remark}

The edge-replacement operation of~\cref{prop:path-replacement} also preserves the $\{a,b\}$-property, even when the replaced edge is simple, as the next lemma shows.

\begin{lemma}\label{lem:leaf-extension}
    Let $a,b$ be distinct real weights.
    Let $K^+$ have a proper $\{a,b\}$-edge-weighting $\omega$, and let $e\in E(K^+)$ be a distinguished edge in the bundle $uv$.
    Put $K=K^+-e$ and let $P$ be an alternately simple path of length $4k+1$ for some positive integer $k$.
    Attach $P$ to $K$ by identifying the two endpoints of $P$ with $u$ and $v$, respectively, and making no other identifications to obtain the graph $K'$.
    Then the restriction of $\omega$ to $K$ extends to a proper $\{a,b\}$-edge-weighting of $K'$ under which every vertex of $K$ has the same weighted degree as in $K^+$.
\end{lemma}

\begin{proof}
    Let $\lambda \in \{a,b\}$ denote the weight of the deleted edge $e$.
    Give the two red edges of $P$ incident with $u$ and $v$ weight $\lambda$.
    These two weights replace the contribution of $e$ at both endpoints.
    Hence, in the resulting weighting of $K'$, every vertex of $K$ has the same weighted degree as in $K^+$, so every edge of $K$ is proper.

    The path $P$ has $2k+1$ red bundles and $2k$ green bundles.
    Alternate the red-edge weights in their order from $u$ to $v$.
    There are $2k$ transitions, so the two prescribed endpoint red weights are both $\lambda$.
    For each green bundle of $P$, the two adjacent red edges have different weights, so the bundle is proper.

    It remains to choose the totals on the green bundles of $P$ so that no red edge of $P$ is conflicting, including the two incident with $u$ and $v$.
    For such a green bundle of multiplicity $m$, the available totals are $\{ma+j(b-a):j\in[0,m]\}$.
    Delete $s_\omega(u)-\lambda$ from the first list and $s_\omega(v)-\lambda$ from the last list when those values are present.
    The endpoint lists remain non-empty, while every internal list has size at least two.
    If one of these green-bundle lists had size at least three before the endpoint deletions, then after the deletions either an endpoint list has size at least two or an internal list has size at least three, so~\cref{lem:list-path}\cref{lem:list-path:large-list} applies.
    Otherwise, all these green bundles are simple, so every list before the endpoint deletions is $\{a,b\}$, and their number $2k$ is even.
    If the two endpoint lists were the same singleton, then the same value would have been deleted from each, giving $s_\omega(u)=s_\omega(v)$ and contradicting properness across the non-empty bundle $uv$ in $K^+$.
    Hence, one can choose distinct values, one from each endpoint list, so~\cref{lem:list-path}\cref{lem:list-path:alternating} applies.
    The chosen consecutive totals make every red edge of $P$ proper, completing a proper weighting.
\end{proof}

The next two lemmas concern modifications considered by Bensmail, Mc Inerney, and Lyngsie~\cite[Lemma~4]{bensmail2022oddab}.
Their proofs differ in the construction of the initial graph but use the same attachment argument, which we describe first.

Root the incidence tree at either a cycle node or a green-bundle node.
Use the augmented multiplicity $m+t$ from the proof of~\cref{prop:path-replacement} for every green-bundle node other than the root.
Once the initial graph has been properly weighted, add the remaining cycle subgraphs in nondecreasing order of the distances of their nodes from the root by the edge replacements described in that proof.
The initial graph is the modified root cycle subgraph or, for a green-bundle root, the modified union of its incident cycle subgraphs.
Thus every remaining attachment uses an extra edge of a non-root bundle.
By~\cref{lem:leaf-extension}, each attachment extends the proper weighting while preserving the weighted degrees of all previously introduced vertices.

\begin{lemma}\label{lem:red-modification}
Let $a,b$ be distinct real weights.
Let $G_0$ be an odd multi-cactus with a specified pasted-cycle representation, and assume that $G_0$ is not a simple cycle.
Let $e$ be a red edge-bundle in this representation.
Increase the multiplicity of $e$ to any value at least two.
The resulting graph has a proper $\{a,b\}$-edge-weighting.
\end{lemma}

\begin{proof}
Root the incidence tree at the cycle node whose cycle subgraph contains $e$.
Take its cycle subgraph as the initial graph, increase the multiplicity of $e$ as prescribed, and use the augmented multiplicities for its green bundles.
If there is only one cycle node, then some green bundle has multiplicity at least two, since $G_0$ is not a simple cycle and all red bundles are simple.
Otherwise, at least one green-bundle node adjacent to the root has a child cycle node, so its augmented multiplicity is at least two.
In either case, both alternating classes of the initial graph contain a multiple bundle, since the multiplicity of $e$ has been increased.
Hence~\cref{lem:cycle-bundle} gives a proper weighting of the initial graph, and the attachment construction above gives a proper weighting of the desired graph.
\end{proof}

\begin{lemma}\label{lem:green-modification}
Let $a,b$ be distinct real weights.
Let $G_0$ be an odd multi-cactus with a valid pasted-cycle representation, and let $uv$ be an edge-bundle of multiplicity one represented by a pasted green-bundle node in that representation.
Replace its unique edge by a path of length five.
The resulting graph has a proper $\{a,b\}$-edge-weighting.
\end{lemma}

\begin{proof}
Relabel and negate the weights if necessary so that $a < b$ and $0 < b$ hold.
Root the incidence tree of the given representation at the green-bundle node $g$ representing $uv$.
Let $r\geq2$ be the number of cycle nodes incident with $g$.
Their cycle subgraphs intersect pairwise precisely in the bundle $uv$, since two cycle nodes cannot have two common neighbors in the incidence tree.
The initial graph is obtained from their union by replacing the common bundle with the prescribed path and giving all other green bundles their augmented multiplicities.

Write the new path as $ux_1y_1x_2y_2v$ and give its edges, in order from $u$ to $v$, the weights $b,b,a,a,b$.
For each cycle node incident with $g$, alternate the red-edge weights along the alternately simple path obtained from its cycle subgraph by deleting the bundle $uv$, starting with weight $b$.
Since this path contains an odd number of red edges, both red edges incident with $u$ and $v$ receive weight $b$.
Then $s(u)=s(v)=(r+1)b$, while the weighted degrees at $x_1,y_1,x_2,y_2$ are, in order, $2b,a+b,2a,a+b$.
Since $b>0$, $a<b$, and $r\geq2$, we have $(r+1)b>2b>a+b>2a$, so every edge of the replacement path is proper.

For each of these cycle nodes independently, choose the totals on the green bundles of its cycle subgraph other than $uv$.
The lists for these green bundles, ordered along the alternately simple path from $u$ to $v$, are even in number.
Delete the forbidden total $rb=(r+1)b-b$ from each endpoint list when it is present.
If one of these lists had size at least three before the endpoint deletions, then~\cref{lem:list-path}\cref{lem:list-path:large-list} applies.
Otherwise, every green bundle under consideration is simple, so every list is $\{a,b\}$, and $rb>b$ shows that no endpoint value is deleted.
In the latter case, choose distinct values from the two endpoint lists and apply~\cref{lem:list-path}\cref{lem:list-path:alternating}.
The chosen totals make every red edge of these cycle nodes proper, and the alternating red weights make every edge in each of their green bundles proper.
Thus the initial graph is properly weighted, and the attachment construction above gives a proper weighting of the desired graph.
\end{proof}

\section{Reduction steps}\label{sec:redsteps}
    In this section, assume that $a, b$ are distinct, non-zero real numbers with $a+b \neq 0$.
    Throughout this section, we use the usual notation $\delta = |b - a|$.
    We give several types of local configurations in which we can construct a partial edge-weighting and then extend it to a proper $\{a,b\}$-edge-weighting of the entire graph.

\subsection{Suspended paths and cycles of even length}

A \emph{suspended path} in a graph $G$ is a path of length at least two whose internal vertices have degree two in $G$ and whose two distinct endpoints have degree at least three.
A \emph{suspended cycle} is a cycle with exactly one vertex of degree at least three in $G$, all other vertices of the cycle having degree two in $G$.
We state two separate lemmas in the same-sign and opposite-sign cases, respectively.

In the same-sign case, we establish a stronger lemma that, in fact, admits a more straightforward proof due to the uniform sign of the weights.

\begin{lemma}\label{lem:samesign-evensuspended}
    Let $a,b$ be distinct non-zero real weights of the same sign.
    Every connected bipartite multigraph containing a suspended path of length two or four with at most one endpoint being a cutvertex has the $\{a,b\}$-property.
    Every connected bipartite multigraph containing a suspended cycle of length four has the $\{a,b\}$-property.
\end{lemma}

\begin{proof}
    Let $G$ be a connected bipartite graph.
    We may assume that $G$ has even order, since otherwise it has the $\{a,b\}$-property by~\cref{cor:odd-order}.
    We may also assume that $0 < a < b$.
    Let $H$ be the graph obtained from $G$ by deleting the internal vertices of the suspended path or the three degree-two vertices of the suspended cycle.
    In the cycle case, $H$ is connected because the suspended cycle forms an endblock of $G$.
    In the path case, if $H$ were disconnected, then both endpoints of the suspended path would be cutvertices of $G$, contrary to the hypothesis.
    Since the number of deleted vertices is odd, $H$ has odd order.

    In the path case, let $u_1,u_2$ be the endpoints of the suspended path, and in the cycle case, let $u_1=u_2$ be the unique vertex of degree at least $3$ on the suspended cycle.
    Use weights $(a,a)$ or $(b,a,a,b)$ when the path or cycle has length $2$ or $4$, respectively, starting at $u_1$.
    Extend this partial edge-weighting to $G$ by applying~\cref{cor:completion-applied} to $G,H$, and the initial weighting defined on the path or cycle.
    Since the weighted degrees of adjacent internal vertices differ by construction, the only conflicting edges may be the edges incident with $u_1$ or $u_2$ that lie on the path or cycle.
    Such an edge is a conflict if the contribution of the other edges incident with $u_i$ is $a$ in total.
    However, since the weights are positive and $a$ is the smaller, this cannot occur, as $u_i$ has $G$-degree at least $3$.
\end{proof}

For the opposite-sign case, we first prove a lemma that gives integer assignments satisfying prescribed parity conditions while excluding one specified value for each of two total counts.
It also identifies the only possible exception to this conclusion.

\begin{lemma}\label{lem:signed-twoholes}
    Let $t$ be a positive integer.
    For $i\in[1,t]$, let $m_i,n_i\geq1$ and $\tau_i\in\{0,1\}$ be integers, and put $M=\sum_i m_i$, $N=\sum_i n_i$.
    Suppose $K,L$ are integers with
    \[
      2\leq K\leq M/2,\qquad 2\leq L\leq N/2,\qquad
      \sum_i\tau_i\equiv K+L\pmod2.
    \]
    For any integers $x_i,y_i$, write $X=\sum_i x_i$ and $Y=\sum_i y_i$.
    There are integers $x_i,y_i$ satisfying $x_i\in[0,m_i]$, $y_i\in[0,n_i]$, and $x_i+y_i\equiv\tau_i\pmod2$, as well as
    \[
      X\equiv K\pmod2,\qquad Y\equiv L\pmod2,
      \qquad X\neq K,\qquad Y\neq L,
    \]
    unless $M=N=4$ and $K=L=2$.
    If no such integers exist, there are nevertheless integers $x_i,y_i$ satisfying $x_i\in[0,m_i]$, $y_i\in[0,n_i]$, and $x_i+y_i\equiv\tau_i\pmod2$ for every $i$, with $X=Y=3$.
    In this exceptional case, $K=L=2$, so these totals have parity opposite to the prescribed total parities.
\end{lemma}

\begin{proof}
    Suppose no integers $x_i,y_i$ satisfy $x_i\in[0,m_i]$, $y_i\in[0,n_i]$, and $x_i+y_i\equiv\tau_i\pmod2$ for every $i$, together with $X\equiv K\pmod2$, $Y\equiv L\pmod2$, $X\neq K$, and $Y\neq L$.
    Because $m_i,n_i\geq1$, every choice of $x_i\in[0,m_i]$ admits a $y_i\in[0,n_i]$ satisfying $x_i+y_i\equiv\tau_i\pmod2$, and the same holds with the two coordinates exchanged.
    Consequently, every integer in $[0,M]$ can occur as $X$ while all the individual congruences hold.
    The hypothesis $\sum_i\tau_i\equiv K+L\pmod2$ means that $X\equiv K\pmod2$ automatically gives $Y\equiv L\pmod2$.
    In particular, $[0,M]$ contains at least two integers of parity $K$.

    If some $n_i\geq3$, choose $X\neq K$ of parity $K$; the two choices for $y_i$ of its required parity cannot both give $Y=L$.
    Thus $n_i\leq2$ for all $i$, and symmetrically $m_i\leq2$; this implies $t \geq 2$, as $M,N \geq 4$.
    If we have a pair $(m_i,n_i)=(1,2)$, take $x_i=\tau_i$, making $y_i=0,2$ both admissible.
    The sum of the other $x$-coordinates ranges over $[0,M-1]$, so $X$ ranges over $[x_i,x_i+M-1]$, which contains at least two values of parity $K$ since $M\geq4$.
    At least one of them gives $X\neq K$; then the value of $y_i$ can be chosen to avoid $L$.
    The case $(2,1)$ is symmetric.
    When $(m_i,n_i,\tau_i)=(2,2,0)$, the other $x$-coordinates can have either parity because their combined upper bound is $M-2\geq2$.
    Choose their total so that $X\equiv K\pmod2$; then choose $x_i,y_i$ independently from $\{0,2\}$ to avoid both forbidden totals.

    Thus every remaining triple $(m_i,n_i,\tau_i)$ is of the form $(1,1,\tau_i)$ or $(2,2,1)$, so $m_i=n_i$ and $M=N$.
    If a triple $(2,2,1)$ occurs, fix $x_i=1$ and use $y_i\in\{0,2\}$.
    The remaining $x$-coordinates attain every integer in $[0,M-2]$; the parity of their sum must be $K-1$.
    If any such total differed from $K-1$, the two available choices for $y_i$ would then avoid $L$ as well.
    Failure therefore requires $K-1$ to be the only integer of its parity in $[0,M-2]$.
    Since $M\geq4$ and $2\leq K\leq M/2$, this forces $M=4$, $K=2$.
    Interchanging the two coordinates gives $N=4$ and $L=2$.
    The total parity then forces exactly two indices with $\tau_i=1$, since $t \leq 3$ holds in this case.

    Otherwise, all multiplicities are one; thus $t \geq 4$.
    Let $z$ be the number of indices with $\tau_i=1$.
    Here $y_i=1-x_i$ at those indices and $y_i=x_i$ at the others.
    For $z=0$ or $z=M$, taking $X$ to be $0$ or $1$ according to its prescribed parity avoids both totals, since $K,L \in [2,M-2]$.
    If $0<z<M$ and $K$ is odd, then $K\geq3$ and $X=1$ can be attained with both $Y=z-1,z+1$; one of the alternatives avoids $L$, a contradiction.
    Therefore, $K$ is even, and $X=0$ forces $z=L$.
    There are at least two indices with $\tau_i=0$ and at least two with $\tau_i=1$, as $2 \leq L \leq N/2 \leq N-2$.
    If $M$ is odd, setting $x_i=0$ at one index with $\tau_i=0$ or at one index with $\tau_i=1$ and setting all other $x_j=1$ gives $X=M-1>K$ and two different possible values of $Y$.
    Hence, $M$ is even; $X=M$ forces $M-z=L=z$.
    Setting $x_i=1$ at two indices with $\tau_i=0$ and $x_i=0$ elsewhere gives $X=2$, $Y=z+2$, so $K=2$.
    Setting $x_i=0$ at two indices with $\tau_i=1$ and $x_i=1$ elsewhere gives $X=M-2$, $Y=z+2$, so $M-2=K=2$.
    Again, $M=N=4$, $K=L=2$, and exactly two indices have $\tau_i=1$.

    In either exceptional case, assign $(1,0)$ and $(0,1)$ at the two indices with $\tau_i = 1$ and $(0,0)$ elsewhere.
    The totals are $(1,1)$.
    Replacing each pair by $(m_i-x_i,n_i-y_i)$ preserves the congruence $x_i+y_i\equiv\tau_i\pmod2$ because $m_i+n_i$ is even for every remaining triple.
    The total of each coordinate changes from $1$ to $M-1=N-1=3$, giving $(X,Y)=(3,3)$.
\end{proof}

We now turn to the opposite-sign case, using~\cref{lem:signed-twoholes} to handle the possible conflicts at the two endpoints of a suspended path.

\begin{lemma}\label{lem:opsign-evensuspended}
    Let $a,b$ be non-zero real weights of opposite sign with $a+b \neq 0$.
    Every connected bipartite multigraph containing a suspended path of length two or four whose endpoints are not cutvertices has the $\{a,b\}$-property.
    Every connected bipartite multigraph containing a suspended cycle of length four has the $\{a,b\}$-property.
\end{lemma}

\begin{proof}
    By~\cref{cor:odd-order} and~\cref{lem:irrational-ratio}, respectively, we may assume that the graph $G$ has even order and $a,b$ have a rational ratio.
    Scale and relabel the weights so that $a,b$ are relatively prime integers with $a > 0$ and $|a| < |b|$.
    Since the weights have opposite signs, we have $\delta = a-b=|a| + |b|$.

    First, we prove the case where we have a suspended path.
    Let $u,v$ be the endpoints, each of degree at least three.
    Delete the one or three internal vertices of the path, obtaining a connected odd-order graph $H$.
    Both $H-u$ and $H-v$ are connected: after deleting an endpoint in the original graph, the deletion of the remaining internal vertices preserves connectivity.
    The vertices $u,v$ are non-adjacent in $H$ because they lie in the same bipartition class.

    For length two, preweight the path $(b,b)$; for length four, use $(a,b,b,a)$.
    Apply~\cref{cor:completion-applied} to $G,H$, and the initial weighting of the path.
    Retain the prescribed parities of the $b$-degrees in $H$ obtained in the proof of~\cref{lem:completion-general}.
    For any weighting $\omega_1$ of $H$ satisfying these parity prescriptions, the weighted-degree difference across each edge induced by $V(H)$ in the resulting weighting of $G$ is not an integral multiple of $2\delta$.
    By construction, an edge connecting two internal vertices of the path cannot be conflicting.
    It remains to ensure that the path-edges incident with one of the endpoints are proper.
    For an endpoint $z\in\{u,v\}$, a conflict is possible only if the total contribution of the edges incident with $z$ in $H$ is exactly $b$; that is, $s_{\omega_1}(z) = ad_H(z) - \delta \db_{\omega_1} (z) = b = a- \delta$.
    Since $a$ and $\delta$ are relatively prime, this can only happen if there is an integer $\ell_z\geq1$ such that
    \[
      d_H(z)=1+\delta\ell_z,\qquad \db_{\omega_1}(z)=1+a\ell_z.
    \]
    By hypothesis on the weights, we have $\delta \geq 2a + 1$, so $2\leq 1 + a\ell_z\leq(1+\delta \ell_z)/2$.
    If $d_H(z) \not \equiv 1 \pmod{\delta}$ or the prescribed parity differs from the parity of the forbidden $b$-degree $1+a \ell_z$, then the weighted degree at the endpoint cannot equal that at the adjacent internal vertex.
    If neither endpoint can attain its forbidden $b$-degree, then the initial completion is already proper.
    If only one endpoint can attain its forbidden $b$-degree, then rename the endpoints if necessary so that this is $u$.
    Since $H-u$ is connected, any number of edges incident with $u$ of the prescribed parity can receive weight $b$, with the remaining prescribed parities realized inside $H-u$.
    To see this, change the prescriptions on the other endpoints of the chosen $b$-edges and apply~\cref{lem:parity-subgraph} to $H-u$.
    Since the original prescription was feasible and the number of chosen $b$-edges satisfies the prescribed parity, we obtain a feasible prescription on $H-u$.
    The forbidden $b$-degree lies between two and half the available degree, so another available $b$-degree of the same parity exists (e.g., $0$ or $1$).
    Choosing such a $b$-degree and realizing the remaining parity prescriptions gives a proper weighting of $G$ in this case.

    If both endpoints can attain their forbidden $b$-degrees, then denote these $b$-degrees by $K=a\ell_u+1$ and $L=a\ell_v+1$; their prescribed parities agree with these $b$-degrees.
    Every component $C_i$ of $H-u-v$ meets both endpoints, since $H-u$ and $H-v$ are connected.
    Let $m_i,n_i\geq1$ count its edges to $u,v$, respectively, and let $\tau_i$ be the sum of its prescribed parities modulo two.
    Note that $d_H(u) = \sum_i m_i = M$ and $d_H(v) = \sum_i n_i = N$.
    Since the parity prescription is feasible on $H$, its total is even, giving $\sum_i\tau_i\equiv K+L\pmod2$.
    Let $x_i$ and $y_i$ denote the number of $b$-edges from $u$ to $C_i$ and from $v$ to $C_i$, respectively.
    Write $X=\sum_i x_i$ and $Y=\sum_i y_i$.
    By~\cref{lem:parity-subgraph}, these extend inside $C_i$ while attaining the parity prescriptions exactly when $x_i+y_i\equiv\tau_i\pmod2$.
    The preceding degree bounds give $2\leq K\leq\sum_i m_i/2$ and $2\leq L\leq\sum_i n_i/2$.
    Therefore,~\cref{lem:signed-twoholes} applies.
    If $x_i,y_i$ can be chosen with the prescribed parities and $X\neq K$, $Y\neq L$, then use them with the original path pattern to obtain a proper weighting.
    Otherwise,~\cref{lem:signed-twoholes} forces $M=N=4$ and $K=L=2$, and supplies $X=Y=3$.
    Since $M=d_H(u)=1+\delta \ell_u$ and $N=d_H(v)=1+ \delta \ell_v$, these equalities force $\delta=3$, $a=1$, and $\ell_u=\ell_v=1$; in particular, both endpoints have degree five in the original graph $G$.

    To use these exceptional boundary $b$-degrees, change the length-two pattern to $(a,a)$ or the length-four pattern to $(b,a,a,b)$.
    Toggle the prescribed parity at both $u$ and $v$ and leave it unchanged elsewhere.
    The total prescribed parity remains even because the two endpoint parities are both reversed.
    At each endpoint, its fixed contribution changes by $\pm \delta$, while the prescribed parity of its number of $b$-incidences also changes, so the parity conditions ensuring unequal weighted degrees within $H$ remain valid.
    The internal component parities $\tau_i$ are unchanged.
    The new forbidden $b$-degree is $1$ at each endpoint, as that is the only $b$-degree inducing a total contribution $a$ in $H$.
    The boundary $b$-degrees $X=Y=3$ from~\cref{lem:signed-twoholes} differ from the new forbidden $b$-degree $1$ at both endpoints.
    All internal edges of either path pattern are proper.
    This completes the suspended path case.

    In the suspended cycle case, let $u$ be the unique vertex of degree at least three on the cycle and let $H$ be the connected graph obtained by deleting the other three vertices.
    First preweight the cycle with $(a,b,b,a)$ in cyclic order starting at $u$.
    Apply~\cref{cor:completion-applied} to $G$, the odd-order subgraph $H$, and the initial weighting of the cycle.
    The two cycle edges at $u$ can be conflicting only if the remainder contributes $b-a=-\delta$ at $u$.
    Changing the cycle pattern to $(b,a,a,b)$ decreases the fixed contribution at $u$ by $2\delta$.
    Every weighted-degree difference across an edge of $H$ therefore changes by an even integer multiple of $\delta$ and remains non-zero by~\cref{cor:completion-applied}, while a cycle edge can be conflicting only if the contribution from $H$ at $u$ is $a-b=\delta$.
    Both conflicts cannot occur, as the contribution from $H$ is unchanged.
\end{proof}

\subsection{A generalized theta graph as an endblock}

A graph is called a \emph{generalized theta graph} if it consists of two designated vertices connected by $k \geq 3$ internally vertex-disjoint paths.
If each of these paths has length $1$ or $3$, we say that we have a \emph{generalized $1$-$3$-theta graph}.
Note that generalized $1$-$3$-theta graphs are bipartite and the two designated vertices are in different bipartition classes.

\begin{lemma}\label{lem:theta-endblock}
  Let $G$ be a connected bipartite multigraph that contains a generalized $1$-$3$-theta graph as an endblock.
  Suppose that at least one of the paths joining the two designated vertices has length $3$.
  If this endblock contains a cutvertex, assume that the cutvertex is one of the designated vertices.
  Then $G$ has the $\{a,b\}$-property for any distinct non-zero $a,b$ with $a+b \neq 0$.
\end{lemma}

\begin{proof}
    By~\cref{cor:odd-order}, we may assume that $G$ has even order.
    We may also assume that $a > 0$ and $|a| < |b|$ both hold: if necessary, relabel and negate the weights to ensure this.
    Label the designated vertices $u$ and $v$ so that if one of them is a cutvertex, then it is $u$.
    Note that we have $d_G(v) \leq d_G(u)$.
    Let $ux_iy_iv$ denote the paths of length $3$, where $i\in[1,r]$, and let $m$ denote the number of edges joining $u$ and $v$.
    By the assumption of the lemma, we have $r \geq 1$, whereas $m = 0$ is allowed.

    Use initial weight $a$ on each edge incident with $v$.
    Since $G$ has even order and $v$ is not a cutvertex, deleting $v$ leaves a connected odd-order graph $H$.
    Apply~\Cref{cor:completion-applied} to $G, H$, and the initial weighting outside $H$.
    With respect to the resulting edge-weighting of $G$, the only conflicting edges may be the ones incident with $v$.
    Since each internal edge $x_iy_i$ is proper and $y_iv$ has weight $a$, the edge $ux_i$ has weight $b$ for every $i\in[1,r]$.

    If $0 < a < b$, then switch the weight of the edges of the form $x_iy_i$ if necessary so that each of them has weight $a$.
    We show that the obtained weighting of $G$ is proper.
    Only the edges of the form $x_i y_i$ were switched.
    Thus, only the weighted degrees of the vertices $x_i$ and $y_i$ may have changed, and, for a fixed $i$, they changed by the same value.
    Therefore, the edges $x_i y_i$ and the edges outside the endblock (if the endblock is not $G$ itself) remain proper.
    The edges $y_i v$ cannot be conflicting as $s(y_i) = 2a < a d_G(v) = s(v)$.
    Similarly, since $a < b$, we have $s(x_i) = a + b < a (d_G(u) -1)+ b \leq s(u)$ for each $i\in[1,r]$.
    The remaining edges are those that connect $u$ and $v$.
    Since $r \geq 1$, the vertex $u$ has at least one incident edge with weight $b$.
    Hence $s(v) = a d_G(v) \leq ad_G(u) < s(u)$.

    If $b < 0 < a$, then the edges of the form $y_i v$ cannot be conflicting, as $s(y_i) \leq 2a < a d_G(v) = s(v)$.
    We may therefore assume that $m\geq1$ and $s(u)=s(v)$, since otherwise the weighting is proper.
    Switch the weights on each path of length $3$ joining $u$ and $v$, and switch the weight of each edge of the bundle $uv$ except possibly one so that the total number of switched edges incident with $v$ is even.
    The switched edges form an even-degree subgraph and thus decompose into edge-disjoint cycles, so by~\cref{cor:completion-applied}, the only conflicting edges still must be incident with $v$.

    In the weighting obtained by the completion, each path of length $3$ had weights $(a,\lambda,b)$ in order from $v$ to $u$ for some $\lambda \in \{a,b\}$.
    Therefore, switching each such path decreases $s(v)$ by $\delta$ and increases $s(u)$ by $\delta$.
    Switching edges joining $u$ and $v$ preserves the difference of the weighted degrees of the endpoints.
    Therefore, $uv$ cannot remain conflicting after switching, as $r \geq 1$.
    After the switches, each $s(y_i)$ belongs to $\{2b,a+b\}$, and all edges incident with $v$ except possibly one have weight $b$.
    If $d_G(v)\geq4$, then $s(v)\leq a+3b<2b\leq s(y_i)$, since $a+b<0$.
    If $d_G(v)=3$, then exactly two edges incident with $v$ were switched, giving $2b<a+2b=s(v)<a+b$.
    Hence none of the edges $y_i v$ is conflicting, completing the proof.
\end{proof}

\subsection{Two-vertex endblocks}

We call a vertex \emph{pendant} if it has only one neighbor.
The first lemma shows that a connected bipartite multigraph with at least three vertices and a pendant vertex admits a proper $\{a,b\}$-edge-weighting when $a$ and $b$ are distinct and have the same sign.
The second lemma applies to distinct non-zero $a,b$ with $a+b\neq0$, provided some vertex $v$ has a pendant neighbor and $G-v$ has at most one component with more than one vertex.

\begin{lemma}\label{lem:pos-nontriv-endblock}
  Let $G$ be a connected bipartite multigraph with at least $3$ vertices that contains a pendant vertex or, equivalently, a two-vertex endblock.
  Then $G$ has the $\{a,b\}$-property for any distinct non-zero $a,b$ of the same sign.
\end{lemma}

\begin{proof}
    Without loss of generality, we may assume $0 < a < b$.
    Also assume that $G$ has even order, because otherwise the statement follows from~\cref{cor:odd-order}.
    Let $u$ be a vertex with a unique neighbor $v$.
    Fix weight $a$ on each edge incident with $u$ and write $H = G-u$.
    Since $G$ has at least $3$ vertices, $v$ has a strictly larger $G$-degree than $u$.
    Applying~\cref{cor:completion-applied} to $G$, the odd-order connected subgraph $H$, and the initial weighting defined on the bundle $uv$ provides an $\{a,b\}$-edge-weighting $\omega$ of $G$ with every edge not incident with $u$ proper.
    The bundle $uv$ is also proper as $s_\omega(v) \geq ad_G(v) > ad_G(u) = s_\omega(u)$ holds; hence, $\omega$ is proper and the proof is complete.
\end{proof}

\begin{lemma}\label{lem:nontriv-endblock}
    Let $G$ be a connected bipartite multigraph with at least $3$ vertices.
    Suppose $G$ has a vertex $v$ such that at least one component of $G-v$ consists of a single pendant vertex and at most one component contains more than one vertex.
    Then $G$ admits a proper $\{a,b\}$-edge-weighting for any distinct non-zero $a,b$ with $a+b \neq 0$.
\end{lemma}

\begin{proof}
    Suppose, for a contradiction, that $G$ admits no proper weighting.
    Then $G$ has even order, and $a,b$ have opposite sign by~\cref{cor:odd-order} and~\cref{lem:pos-nontriv-endblock}, respectively.
    By~\cref{lem:irrational-ratio}, we may assume that the weights have a rational ratio; furthermore, they are relatively prime integers.
    Relabel and negate the weights if necessary so that $b < 0 < a$ and $|a| < |b|$ hold.
    If each component of $G-v$ is a single vertex, then the underlying simple graph of $G$ is a star with at least $3$ vertices.
    Assigning weight $a$ to each edge gives a proper $\{a,b\}$-edge-weighting.
    Hence, $G-v$ has a unique component with at least two vertices.
    Let $u_1,\dots,u_r$ be all pendant neighbors of $v$.

    First, assume that there exists a vertex $u_i$ such that $d_G(u_i) \not\equiv d_G(v) \pmod{\delta}$.
    Then use weight $a$ on each edge incident with $u_i$ and extend this weighting using~\cref{cor:completion-applied} on $G$, the connected odd-order subgraph $G-u_i$, and the initial weighting defined on the bundle $u_iv$.
    The resulting weighting is proper as the necessary condition of a conflict given by~\cref{lem:conflict-necessary} fails on the bundle $u_i v$, yielding a contradiction.
    Henceforth, we assume that the degrees of $u_1,\dots,u_r$ and $v$ are in the same residue class modulo~$\delta$.

    Relabel the pendant vertices if necessary so that $u_1$ is a minimum degree pendant vertex.
    Use weight $a$ on exactly one edge incident with $u_1$ and weight $b$ on every other edge incident with $u_1$.
    Extend this weighting to $G$ using~\cref{cor:completion-applied} on $G$, the connected subgraph $G-u_1$ of odd order, and the initial weighting defined on the bundle $u_1v$.
    We may assume that $u_1 v$ is a conflict with respect to the resulting weighting $\omega$, since otherwise $\omega$ is proper.
    Put $H = G-\{u_1,\dots,u_r\}$.
    Then $H$ has at least three vertices and $v$ is not a cutvertex in $H$.
    It follows that for any two different edges $e_1, e_2$ of $H$ incident with $v$, there exists a cycle in $H$ containing both of them.
    By~\cref{cor:completion-applied}, switching any cycle in $G$ preserves that the edges of $G-u_1$ are proper.

    If there exist two edges $e_1,e_2$ of the same weight incident with $v$ in $H$, then switching a cycle containing both edges results in a proper weighting, as $s_\omega(u_1)$ remains unchanged while $s_\omega(v)$ changes, a contradiction.
    It follows that $d_H(v) \leq 2$.
    If $d_H(v) = 2$, then the two edges incident with $v$ in $H$ have different weights.
    Hence, we have $r \geq 2$, since otherwise $d_H(v) = d_G(v) - d_G(u_1) \geq  \delta \geq 3$ would hold.
    Relabel the pendant vertices other than $u_1$ so that $u_r$ is a maximum degree pendant vertex.
    If $u_r$ has degree at least $3$, then we can obtain a proper weighting by switching two edges of the bundle $u_r v$ of the same weight; this cycle switch preserves the weighted degree of $u_1$ while changing the weighted degree of $v$.
    Hence, each $u_i$ has degree at most $2$.
    If $u_r$ has degree $2$, then each pendant vertex has degree $2$, as their degrees belong to the same residue class modulo~$\delta$.
    Since the bundle $u_1 v$ has one edge of each admissible weight and $u_r v$ is not conflicting, the bundle $u_r v$ consists of two edges of the same weight.
    Similarly, switch these two edges to obtain a proper weighting of $G$ and reach a contradiction.

    The remaining case is where $u_r$, and therefore each pendant vertex, has degree $1$.
    Note that $u_1v$ has weight $a$ and $u_i v$ has weight $b$ for each $i\in[2,r]$ as none of these edges is conflicting.
    Then the total contribution of the edges $u_2 v, \dots, u_r v$ to the weighted degree of $v$ is $(r-1)b$, where $r-1 \geq 1$.
    As $u_1 v$ is a conflict by assumption, the sum of edge-weights incident with $v$ in $G-u_1$ is $0$.
    Since $|a| < |b|$ and $b < 0 < a$, there must exist at least two edges incident with $v$ in $H$ with weight $a$ to compensate the contribution $(r-1)b \leq b$ from the pendant edges.
    This is a contradiction, since we observed that there cannot be two edges of the same weight incident with $v$ in $H$; the proof of the lemma is complete.
\end{proof}

\subsection{Local extrema}

Let $G$ be a multigraph.
We call a vertex $v$ a \emph{local minimum (maximum)} among a set of vertices $U \subseteq V(G)$ if each $u \in U$ satisfies $d_G(v) \leq d_G(u)$ ($d_G(v) \geq d_G(u)$).
We jointly refer to a local minimum or maximum as a \emph{local extremum}.
Analogously, $v$ is a \emph{weighted local minimum} or \emph{maximum} if the corresponding inequalities hold for weighted degrees instead of degrees.
We say that $v$ is a \emph{strict} local minimum or maximum if the corresponding inequality holds strictly.
A frequently used idea is to select a vertex $w$ of $G$, define the weights of the edges incident with $w$, and extend this weighting to $G$ by applying~\cref{cor:completion-applied} to $G, H = G-w$, and the initial weighting defined on the edges incident with $w$.
To obtain a proper weighting, it suffices to avoid conflicting edges incident with $w$.
A way to achieve this is to make $w$ a strict weighted local extremum among certain neighbors.
Building on this idea, we can establish a lemma in the same-sign and opposite-sign cases, respectively.

We introduce the notion of \emph{problematic} edge-bundles in both cases.
If $a$ and $b$ have the same sign, then we call an edge-bundle $uv$ problematic if $d_G(u) = d_G(v)$ and there is no suspended path of length $3$ connecting $u$ and $v$.
If $a$ and $b$ have opposite signs, then $uv$ is problematic if $d_G(u) \equiv d_G(v) \pmod{\delta}$ and, if $d_G(u) = d_G(v)$, then there is no suspended path of length $3$ connecting $u$ and $v$.
In the same-sign case, this definition applies to arbitrary distinct real weights of the same sign, whereas in the opposite-sign case, we use it only for relatively prime integer weights.
Introduce the following notation:
\[
A_G(w) \coloneqq
\begin{cases}
    \{v \in N_G(w): vw \text{ is problematic}\}& \text{if } a,b \text{ have the same sign,}\\
    \{v \in N_G(w): vw \text{ is problematic}, d_G(v) \geq d_G(w)\}& \text{if } a,b \text{ have opposite signs.}
\end{cases}
\]
We say that a suspended path is incident with a vertex if that vertex is an endpoint of the path.
Similarly, a suspended cycle is incident with a vertex if that vertex is the unique vertex on the cycle with degree at least $3$.
Let $G^*(w)$ denote the graph we obtain from $G$ by removing $w$ and the internal vertices of suspended paths of length $3$ incident with $w$.

\begin{lemma}\label{lem:samesign-locminmax}
    Let $G$ be a connected bipartite loopless multigraph, and let $a,b$ be distinct real weights of the same sign.
    Let $w \in V(G)$ be a non-cutvertex of $G$-degree at least $3$.
    Suppose that $w$ is a local extremum among its neighbors of $G$-degree at least $3$ and that each suspended path incident with $w$ has length $3$.
    If there exists a set of edges $F \subseteq E(G^*(w))$ such that $G^*(w)-F$ is connected and $d_F(v) \geq 1$ holds for each vertex $v \in A_G(w)$, then $G$ has the $\{a,b\}$-property.
\end{lemma}

\begin{proof}
    Without loss of generality, we may assume that $0 < a < b$.
    By~\cref{cor:odd-order}, it suffices to verify the statement if $G$ has even order.
    Since $w$ is not a cutvertex and has degree at least $3$, it has no neighbor of degree $1$.
    Starting at $w$ and following consecutive edges through vertices of degree $2$ must lead to another vertex of degree at least $3$, since ending at a vertex of degree $1$ or returning to $w$ would make $w$ a cutvertex.
    Consequently, by hypothesis, every neighbor of $w$ with degree $2$ lies on a suspended path of length $3$ incident with $w$.
    Also note that the graph $G^*(w)$ has odd order, as $G$ has even order and each suspended path of length $3$ has $2$ internal vertices, giving an odd number of deleted vertices in total.
    The general idea is to make $w$ a strict weighted local extremum among its neighbors of $G$-degree at least $3$ and ensure that the neighbors of $w$ with degree $2$ have a strictly smaller weighted degree than that of $w$.

    First, assume that $w$ is a local minimum.
    Define initial weight $a$ on each edge incident with $w$ and weights $(a,a,b)$ on each suspended path of length $3$ incident with $w$.
    Use weight $b$ on every edge of $F$.
    Extend this weighting to $G$ by applying~\cref{cor:completion-applied} to $G$, the connected graph $G^*(w)-F$ of odd order, and the initial weighting defined outside $G^*(w)-F$.
    By the corollary, each edge of $G^*(w)$ is proper with respect to the resulting weighting $\omega$, so it suffices to verify that the edges outside $G^*(w)$ cannot be conflicting.
    First, consider a bundle $wv$ with $d_G(v) \geq 3$.
    \begin{itemize}
        \item If $d_G(w) < d_G(v)$, then we have $s_\omega(w) = ad_G(w) < ad_G(v) \leq s_\omega(v)$.
        \item Otherwise $d_G(w) = d_G(v)$.
        If at least one suspended path of length $3$ connects $w$ and $v$, then the edge with endpoint $v$ on that path has weight $b$.
        Otherwise, $wv$ is problematic and therefore $v \in A_G(w)$, so $v$ has an incident edge in $F$ of weight $b$.
        In either case, we have $s_\omega(w) = ad_G(w) = ad_G(v) < s_\omega(v)$.
    \end{itemize}
    Now consider a suspended path $wxyv$ incident with $w$.
    \begin{itemize}
        \item The edge $xy$ is proper, as $wx$ and $yv$ have different weights.
        \item We have $s_\omega(w) = ad_G(w) > 2a = s_\omega(x)$, and similarly $s_\omega(v) \geq a(d_G(v) -1) + b > a+b = s_\omega(y)$.
          Hence, the outer edges are also proper.
    \end{itemize}

    Now assume that $w$ is a local maximum.
    Use weight $b$ on each edge incident with $w$, weights $(b,a,a)$ on each suspended path of length $3$ incident with $w$, and weight $a$ on each edge of $F$.
    Apply~\cref{cor:completion-applied} as in the local-minimum case and denote the resulting weighting by $\omega$.
    For every neighbor $v$ of $w$ with $d_G(v) \geq 3$, we have $s_\omega(v) \leq bd_G(v) \leq bd_G(w)=s_\omega(w)$.
    As in the local-minimum case, if $d_G(v)=d_G(w)$, then $v$ has an incident edge of weight $a$ on a suspended path connecting $w$ and $v$ or in $F$.
    Hence, $s_\omega(v)<s_\omega(w)$.
    For a suspended path $wxyv$ incident with $w$, we have $s_\omega(w)=bd_G(w)>a+b=s_\omega(x)>2a=s_\omega(y)$ and $s_\omega(v)\geq ad_G(v)>2a=s_\omega(y)$.
    Thus, every edge outside $G^*(w)$ is proper, and $\omega$ is a proper weighting of $G$.
\end{proof}

In the opposite-sign case, the preceding argument does not guarantee that a local minimum can be made a strict weighted local minimum.
The following lemma gives a sufficient condition for obtaining a proper weighting by making a vertex a strict weighted local maximum.

\begin{lemma}\label{lem:opsign-locmax}
    Let $G$ be a connected bipartite loopless multigraph, and let $a \neq b$ be non-zero relatively prime integer weights of opposite signs.
    Let $w \in V(G)$ be a non-cutvertex of $G$-degree at least $3$.
    Assume that there exists a set of edges $F \subseteq E(G^*(w))$ such that $G^*(w)-F$ is connected and
    \[
    d_F(v) \geq \frac{|a|}{\delta}(d_G(v)-d_G(w)) +1
    \]
    holds for each vertex $v \in A_G(w)$, where $\delta = |a-b|$.
    Then $G$ has the $\{a,b\}$-property.
\end{lemma}

\begin{proof}
    By~\cref{cor:odd-order}, it suffices to verify the statement for even-order $G$.
    Negate the weights if necessary so that $b < 0 < a$ holds.
    Hence, $|a| = a$ and $\delta = a-b$.
    We use a similar argument as in the proof of~\cref{lem:samesign-locminmax}.
    By~\cref{lem:conflict-necessary}, we only need to focus on the neighbors of $w$ with $G$-degree in the same residue class modulo~$\delta$ as that of $w$.
    We aim to make $w$ a strict weighted local maximum among such neighbors.

    Define initial weight $a$ on each edge incident with $w$ and weight $b$ on each edge of $F$.
    Since $G^*(w)-F$ is connected, $G-w-F$ is also connected.
    Now, we apply~\cref{cor:completion-applied} to $G$, the connected odd-order subgraph $G-w-F$, and the weighting defined on the edges in $F$ and the edges incident with $w$.
    For the resulting weighting $\omega$ of $G$, we have $s_\omega(w)=ad_G(w)$, and the edges not incident with $w$ are proper, so it suffices to show that the edges incident with $w$ cannot be conflicting either.
    Note that the properness of the internal edge $xy$ of a suspended path $wxyv$ of length $3$ implies $\omega(yv) = b$, since $\omega(wx) = a$.
    Let $wv$ be an arbitrary edge-bundle incident with $w$.
    \begin{itemize}
        \item If $d_G(w) > d_G(v)$, then we have $s_\omega(w) = ad_G(w) > ad_G(v) \geq s_\omega(v)$.
        \item If $d_G(w) \not \equiv d_G(v) \pmod{\delta}$, then the necessary condition for $wv$ to be conflicting fails, so $wv$ is proper.
        \item If $d_G(w) = d_G(v)$ and at least one suspended path of length $3$ connects $w$ and $v$, then the edge with endpoint $v$ on that path has weight $b$.
        Thus, we have $s_\omega(w) = ad_G(w) = ad_G(v) > s_\omega(v)$.
        \item Otherwise, $d_G(w) \leq d_G(v)$ and $wv$ is problematic, so $v \in A_G(w)$.
        Let $\alpha_v = d_G(v) -d_G(w) \geq 0$.
        By construction, we have $\db_\omega(v) \geq d_F(v) > \frac{a}{\delta}\alpha_v$.
        Recall that we can express the weighted degree of $v$ as $s_\omega(v) = ad_G(v) + (b-a) \db_\omega(v) = ad_G(v) - \delta \db_\omega(v)$.
        The properness of $wv$ follows, since
        \[
        s_\omega(w) - s_\omega(v) =
        a(d_G(w) - d_G(v)) + \delta  \db_\omega(v) > a (-\alpha_v) + \delta \frac{a}{\delta} \alpha_v = 0.
        \]
    \end{itemize}
\end{proof}

\section{Proof of the two-weight characterization}\label{sec:proof-completion}

It remains to prove the converse implication of~\cref{thm:main}.
By~\cref{rem:mainthm-irratratio}, it suffices to consider weights with a rational ratio.
Furthermore, without loss of generality, we may assume that $a,b$ are relatively prime integers with $a > 0$ and $|a| < |b|$ (relabel and scale if necessary).
Suppose, for a contradiction, that some connected bipartite loopless multigraph not belonging to the odd multi-cactus family has no proper $\{a,b\}$-edge-weighting for these fixed weights $a,b$.
Choose such a counterexample $G$ with as few vertices as possible and, subject to that, with as few edges as possible.
By~\cref{cor:odd-order}, $|V(G)|$ is even, and $|V(G)|\geq4$ since every connected loopless multigraph on two vertices is $K_2^m$.

\subsection{Properties of the minimal counterexample}

First, we establish two properties of $G$ that follow from minimality.
These strongly depend on the modification lemmas in~\cref{subsec:omc-obs-mod}, which give small modifications of odd multi-cacti resulting in graphs with the $\{a,b\}$-property.

\begin{proposition}\label{prop:minimality-prop}
The connected minimal counterexample $G$ satisfies the following properties.
    \begin{enumerate}[label=(M\arabic*),ref=(M\arabic*)]
        \item\label{min:M1} $G$ does not contain a suspended path or cycle of length at least $5$.
        \item\label{min:M2} $G$ does not contain a bundle $xy$ of multiplicity at least $2$ such that $x$ and $y$ each have exactly two distinct neighbors.
        In particular, $G$ does not contain two adjacent vertices $x,y$ of $G$-degree $q+1$ joined by exactly $q \geq 2$ parallel edges.
    \end{enumerate}
\end{proposition}

\begin{proof}
    We prove the two properties in this order using the minimality property of $G$ and the modification lemmas.

    \smallskip
    \noindent\emph{Proof of~\cref{min:M1}.}
    Assume that $G$ contains a suspended path or cycle of length at least $5$ and choose a five-edge subpath $P=y_1x_1y_2x_2y_3x_3$ whose four internal vertices have degree two.
    In the suspended-cycle case, the cycle has even length, so $P$ is a proper subpath.
    Form $G_0$ by deleting the four internal vertices and adding one marked edge $e=y_1x_3$, even if that bundle already exists.
    Since $P$ has odd length, its distinct endpoints lie in opposite bipartition classes.
    Thus $G_0$ is connected, bipartite, and loopless, with $|V(G_0)|=|V(G)|-4$.

    If $G_0$ is not an odd multi-cactus, minimality gives a proper weighting of $G_0$, and~\cref{lem:leaf-extension} gives a proper weighting of $G$ since $P$ is a simple path of length five.
    This contradicts the choice of $G$, so $G_0$ is an odd multi-cactus.

    Let $p_e$ denote the edge-bundle of $G_0$ containing $e$.
    If $p_e$ is multiple, then~\cref{prop:path-replacement} implies that $G$ is an odd multi-cactus, a contradiction.
    Hence $p_e$ is simple.
    Also, $G_0\neq K_2$, since both endpoints of $P$ have degree at least $2$ in $G$.
    Therefore, $G_0$ admits a pasted-cycle representation by~\cref{prop:models}; fix such a representation.

    If $p_e$ is red or private green, replace it in its unique cycle subgraph by the five simple bundles of $P$, colored alternately starting and ending with the color of $p_e$.
    The cycle length increases by four.
    In the incidence tree, delete the node for $p_e$ if it was private green and add private nodes with multiplicity label $1$ for the new green bundles.
    This gives a pasted-cycle representation of $G$, a contradiction.

    It remains that $p_e$ is a pasted green bundle of multiplicity one.
    In this case,~\cref{lem:green-modification} gives a proper weighting of $G$, again contradicting the choice of $G$.
    This proves~\cref{min:M1}.

\smallskip
\noindent\emph{Proof of~\cref{min:M2}.}
Suppose $G$ contains a bundle $xy$ of multiplicity $q \geq 2$ such that $x$ and $y$ each have exactly two distinct neighbors.
Replace the bundle by one marked edge $e=xy$, obtaining $G_0$.
The endpoints are unchanged, so $G_0$ is connected, bipartite, and loopless, and $|E(G_0)|=|E(G)|-(q-1)<|E(G)|$ while its order is unchanged.

If $G_0$ is not an odd multi-cactus, then the minimality of $G$ gives a proper weighting $\omega$ of $G_0$.
For the extension argument, let $\lambda =\omega(e)$, delete the marked edge $e$, and restore the original bundle of $q$ parallel edges.
Its total can be any of the $q+1$ distinct values $qa,(q-1)a + b,\dots,qb$.
Replacing the contribution $\lambda$ by any such total changes the weighted degrees of $x$ and $y$ by the same amount.
Their difference is therefore unchanged, so the weighted degrees still differ across the bundle $xy$.
At each of $x$ and $y$, the unique incident bundle different from the $xy$-bundle forbids at most one total.
Together, these two bundles forbid at most two of the $q+1$ available totals.
Since $q+1\geq3$, an allowed total exists; using such a total gives a proper weighting of $G$, a contradiction.
Hence, $G_0$ must be an odd multi-cactus.
Also, $|V(G_0)| = |V(G)| \geq 3$, so $G_0$ admits a pasted-cycle representation by~\cref{prop:models}.

Fix a valid pasted-cycle representation of $G_0$.
If the marked bundle is green in this representation, then increasing its multiplicity label gives a valid representation of $G$, contradicting the choice of $G$.
Moreover, $G_0$ cannot be a simple cycle, since one of its two alternating colorings makes the marked bundle green, giving the same contradiction.
Thus $G_0$ is not a simple cycle and the marked bundle is red in the chosen representation, so~\cref{lem:red-modification} gives a proper $\{a,b\}$-edge-weighting of $G$, again contradicting the choice of $G$.
\end{proof}
The lemmas in~\cref{sec:redsteps} give sufficient conditions for the $\{a,b\}$-property.
Since $G$ does not admit a proper $\{a,b\}$-edge-weighting, we deduce the following properties.

\begin{proposition}\label{prop:redsteps}
    The connected minimal counterexample $G$ satisfies the following properties.
    \begin{enumerate}[label=(R\arabic*),ref=(R\arabic*)]
    \item\label{red:R1} $G$ contains no suspended cycle of length $4$.
      If $a,b$ have the same sign, then $G$ contains no suspended path of length $2$ or $4$ with a non-cutvertex endpoint.
      If $a,b$ have opposite signs, then $G$ contains no suspended path of length $2$ or $4$ connecting two non-cutvertices.
        \item\label{red:R2} No endblock of $G$ is a generalized $1$-$3$-theta graph with at least one path of length $3$, where the unique cutvertex of the endblock (if it exists) is a designated vertex of the theta graph.
        \item\label{red:R3} $G$ contains an endblock with at least $3$ vertices.
        \item\label{red:R4} The underlying simple graph of $G$ is not a cycle.
\end{enumerate}
\end{proposition}

\begin{proof}
    Properties~\cref{red:R1} and~\cref{red:R2} follow from~\cref{lem:samesign-evensuspended,lem:opsign-evensuspended} and~\cref{lem:theta-endblock}, respectively.
    If the underlying simple graph of $G$ were a cycle, then it would have even length, and by~\cref{lem:cycle-bundle}, $G$ would either have a proper weighting or would be an alternately simple cycle of length $2$ modulo~$4$, hence an odd multi-cactus.
    Both alternatives contradict the choice of $G$, which proves~\cref{red:R4}.

    It remains to verify~\cref{red:R3}.
    If $G$ is $2$-connected, then $G$ itself is an endblock with at least three vertices.
    Otherwise, there exists a cutvertex $v$ in $G$ such that each block that contains it except possibly one is an endblock of $G$.
    To see this, delete the leaves of the block-cutvertex tree of $G$.
    In the resulting tree, every block node still has degree at least two.
    Choose a cutvertex $v$ whose node has degree at most one in this tree.
    If each endblock that contains $v$ had two vertices, then $G$ would have the $\{a,b\}$-property by~\cref{lem:nontriv-endblock}, a contradiction.
    Therefore, one of the endblocks that contain $v$ has at least $3$ vertices, which completes the proof of~\cref{red:R3}.
\end{proof}

\subsection{Simultaneous edge-deletion while preserving connectivity}

In the rest of the proof, we seek to determine a vertex $w_0$ that satisfies the conditions of~\cref{lem:samesign-locminmax} or~\cref{lem:opsign-locmax} in the same-sign or opposite-sign case, respectively.
Let us call such a vertex $w_0$ an \emph{enforceable extremum}.
If $G$ contains an enforceable extremum, then it has the $\{a,b\}$-property by these lemmas.

Note that since $G$ is bipartite, $N_G(w)$ is independent for any vertex $w$.
Therefore, $A_G(w)$ is also independent, as it is a subset of $N_G(w)$.
Both mentioned lemmas require an edge set that satisfies certain lower bounds on the degrees in an independent set of vertices and leaves a connected graph upon deletion.
We deduce a necessary and sufficient condition for this problem from the upper bound specialization of Frank's degree-constrained spanning tree theorem~\cite[Theorem~9.1.16]{Frank2011}.

\begin{theorem}[Frank~\cite{Frank2011}]\label{thm:frank}
Let $J$ be a connected loopless multigraph, let $A\subseteq V(J)$ be independent, and let $g\colon A\to\Z_{\geq0}$ be a function.
There is a spanning tree $T$ of $J$ with $d_T(x)\leq g(x)$ for every $x\in A$ if and only if, for every non-empty $X\subseteq A$,
\[
  \sum_{x\in X}g(x)\geq |X|+\sigma_J(X)-1.
\]
\end{theorem}

The criterion extends from simple graphs to loopless multigraphs by retaining one edge from each parallel class, since this preserves possible spanning tree degree vectors and the vertex sets of components after the deletion of a subset $X$ remain the same.
Applying the theorem with a bound function of the form $g(x)=d_J(x)-h(x) -1$ gives the following criterion.

\begin{corollary}\label{cor:frank-leaving}
Let $J$ be a connected loopless multigraph, let $A\subseteq V(J)$ be independent, and let $h\colon A\to\Z_{\geq0}$ be a function.
There is a set of edges $F\subseteq E(J)$ such that $J-F$ is connected and $d_F(x)\geq h(x)+1$ holds for every $x\in A$ if and only if, for every non-empty $X\subseteq A$,
\begin{equation}\label{eq:frank-leaving}
  \sum_{x\in X}(d_J(x)-h(x))\geq 2|X|+\sigma_J(X)-1.
\end{equation}
\end{corollary}

\begin{proof}
    If there exists a vertex $x \in A$ such that $h(x) + 1 > d_J(x)$, then the desired set of edges clearly does not exist.
    Inequality~\cref{eq:frank-leaving} also fails for $X = \{x\}$, since the left-hand side is non-positive, while the value of the right-hand side is at least $1$.
    Henceforth, assume $h(x) + 1 \leq d_J(x)$ for each $x \in A$.

    The desired set of edges exists exactly when $J$ has a connected spanning subgraph with degree at most $d_J(x)-h(x)-1 \geq 0$ at each $x\in A$.
    This is equivalent to the existence of a spanning tree satisfying the same constraints, since a spanning tree of an appropriate subgraph also satisfies the upper bounds.
    Apply~\cref{thm:frank} with the given non-negative bounds.
    Rearranging the inequality, we obtain exactly~\cref{eq:frank-leaving}.
\end{proof}

\subsection{The set of possible candidates}

The compression of suspended paths of length $3$ allows us to work in a smaller graph while preserving the degrees of the remaining vertices and the connectivity properties needed for the edge-deletion argument.
First, choose an endblock $B$ of $G$ that has at least $3$ vertices, and let $c_0$ be its cutvertex when $B\neq G$.
Such an endblock exists by~\cref{prop:redsteps}\cref{red:R3}.

Extend each edge of $B$ through vertices of degree two in $G$ as far as possible.
A closed chain would be either a suspended cycle of length at least $4$, excluded by~\cref{prop:redsteps}\cref{red:R1} and~\cref{prop:minimality-prop}\cref{min:M1}, or a bundle of multiplicity two with an endpoint of degree two, which forms a two-vertex block and hence cannot lie in $B$, or the whole graph with every degree two, excluded by~\cref{prop:redsteps}\cref{red:R4}.

Thus, each maximal chain is a path with distinct endpoints of degree at least three in $G$, that is, an edge or a suspended path in $G$.
By~\cref{prop:minimality-prop}\cref{min:M1}, the length of such a path is at most $4$.
The interiors of these maximal paths are disjoint, since at a degree-two vertex the continuation is forced.

Since $B$ contains at most one cutvertex of $G$, each suspended path in $B$ has at least one non-cutvertex endpoint, and both endpoints are non-cutvertices if $c_0$ is absent or the path is not incident with $c_0$.
Hence, we can deduce a stronger property regarding the possible lengths of these paths from~\cref{prop:redsteps}\cref{red:R1}.
\begin{itemize}
    \item If $a,b$ have the same sign, then each suspended path in $B$ has length $3$.
    \item If $a,b$ have opposite signs, then each suspended path in $B$ has length $3$ unless $c_0$ is present and the path is incident with $c_0$.
\end{itemize}

Suppress every three-edge suspended path of $B$ into a \emph{blue edge}, retain every other edge as \emph{white}, and call the resulting graphs $R_B$ and $R$ for the endblock and the whole graph, respectively.
Parallel edges, including white and blue edges with the same endpoints, are retained.
The bipartition is preserved because each suppressed path has odd length.
The graph $R_B$ has at least $3$ vertices, since otherwise $B$ would either have at most $2$ vertices, or $B$ would be a suspended cycle, excluded by~\cref{prop:redsteps}\cref{red:R1} and~\cref{prop:minimality-prop}\cref{min:M1}, or a generalized theta graph with at least one path of length $3$, excluded by~\cref{prop:redsteps}\cref{red:R2}.
The graph $R_B$ is loopless and $2$-connected.
Looplessness follows from the fact that the suppressed paths have distinct endpoints.
To verify $2$-connectivity, a path in $B-v$ between vertices retained after compression cannot traverse the interior of a suppressed path ending at $v$, since those internal vertices have no other connection in $B-v$.
Replacing the other suppressed paths on such a path by their blue edges gives a walk in $R_B-v$.
Any part outside $R_B$ remains attached only at $c_0$; therefore, $R_B$ is, in fact, an endblock of $R$, as $c_0$ is a cutvertex in $R$ if it exists.

For every vertex $v \in V(R)$, we have $d_R(v)=d_G(v)$.
We abbreviate this common degree function to $d$.
Every suppression removes two vertices, so $R$ has even order.
If $a,b$ have the same sign, then each vertex in $R_B$ has $R$-degree at least $3$.
If $a,b$ have opposite signs, then vertices of $R$-degree $2$ can occur in $R_B$ only when $c_0$ is present, as internal vertices of suspended paths of length $2$ or $4$ incident with $c_0$.

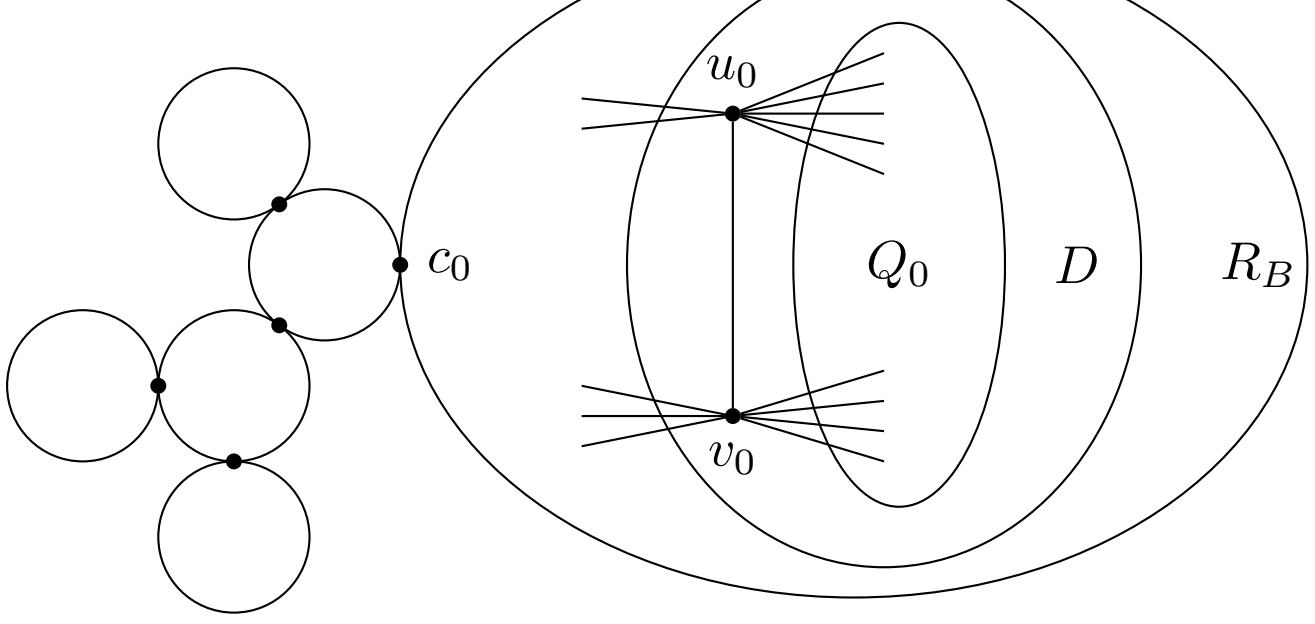
\begin{figure}[H]
    \centering
    \begin{tikzpicture}
        \begin{scope}
            \tikzset{VertexStyle/.style = {shape = circle, fill = black , minimum size = 3pt, inner sep=0pt}}
            \Vertex[x=0,y=1,L=$u_0$, LabelOut = True, Lpos = 90]{u0}
            \Vertex[x=0,y=-1,L=$v_0$, LabelOut = True, Lpos = -90]{v0}
            \Vertex[x=-2.2,y=0,L=$c_0$,LabelOut = True, Lpos = 0]{c0}
            \Vertex[x=-3,y=-0.4,L=$\phantom{}$, LabelOut = True]{r1}
            \Vertex[x=-3,y=0.4,L=$\phantom{}$, LabelOut = True]{r2}
            \Vertex[x=-3.3,y=-1.3,L=$\phantom{}$, LabelOut = True]{r3}
            \Vertex[x=-3.8,y=-0.8,L=$\phantom{}$, LabelOut = True]{r4}

            \coordinate (p1) at (1,1.4);
            \coordinate (p2) at (1,1.2);
            \coordinate (p3) at (1,1);
            \coordinate (p4) at (1,0.8);
            \coordinate (p5) at (1,0.6);
            \coordinate (p6) at (1,-0.7);
            \coordinate (p7) at (1,-0.9);
            \coordinate (p8) at (1,-1.1);
            \coordinate (p9) at (1,-1.3);

            \coordinate (q1) at (-1,1.1);
            \coordinate (q2) at (-1,0.9);
            \coordinate (q3) at (-1,-0.8);
            \coordinate (q4) at (-1,-1);
            \coordinate (q5) at (-1,-1.2);

            \draw[] (u0) -- (v0);
            \draw[] (u0) -- (p1);
            \draw[] (u0) -- (p2);
            \draw[] (u0) -- (p3);
            \draw[] (u0) -- (p4);
            \draw[] (u0) -- (p5);
            \draw[] (v0) -- (p6);
            \draw[] (v0) -- (p7);
            \draw[] (v0) -- (p8);
            \draw[] (v0) -- (p9);
            \draw[] (u0) -- (q1);
            \draw[] (u0) -- (q2);
            \draw[] (v0) -- (q3);
            \draw[] (v0) -- (q4);
            \draw[] (v0) -- (q5);

            \draw[] (1.1,0) ellipse (0.7cm and 1.6cm) node at (1.1,0) {$Q_0$};
            \draw[] (1,0) ellipse (1.7cm and 2cm) node[right] at (2,0) {$D$};
            \draw[] (0.8,0) ellipse (3cm and 2.2 cm) node[right] at (3.1,0) {$R_B$};
            \draw[] (-2.7,0) circle (0.5cm);
            \draw[] (-3.3,0.8) circle (0.5cm);
            \draw[] (-3.3,-0.8) circle (0.5cm);
            \draw[] (-4.3,-0.8) circle (0.5cm);
            \draw[] (-3.3,-1.8) circle (0.5cm);
        \end{scope}
    \end{tikzpicture}
    \caption{Schematic illustration of the case where $c_0$, $u_0$, and $v_0$ are present.}
    \label{fig:candidate-region}
\end{figure}

Next, we define a subgraph of $R_B$ from which to choose a candidate vertex.
We call an edge-bundle $uv$ of $R_B$ problematic in $R$ if $u$ and $v$ are adjacent in $G$ and the bundle $uv$ is problematic in $G$.
\begin{itemize}
    \item If $a,b$ have the same sign, then $uv$ is problematic if $d(u) = d(v)$ and no blue edge connects $u$ and $v$.
    \item If $a,b$ have opposite signs, then $uv$ is problematic if $uv$ contains at least one white edge, $d(u) \equiv d(v) \pmod{\delta}$, and, if $d(u) = d(v)$, then there is no blue edge connecting $u$ and $v$.
\end{itemize}
For a connected graph $H$, let us call an edge-bundle $uv$ \emph{$H$-separating} if $H-u-v$ is disconnected.
If there exists an $R_B$-separating problematic edge-bundle in $R_B$ with endpoints distinct from $c_0$ (if $c_0$ exists), choose $u_0v_0$ to be such a bundle, and a component $Q_0$ of $R_B-u_0-v_0$ not containing $c_0$ (or any component if $c_0$ does not exist) so that $|V(Q_0)|$ is minimum over all such choices.
Put $D=R_B[V(Q_0)\cup\{u_0,v_0\}]$.
Otherwise, put $D=R_B$.
We call the vertices $c_0,u_0,v_0$ that are present the \emph{distinguished vertices}.
In the case with no separating problematic bundle, the only distinguished vertex may be $c_0$.
The case where all three distinguished vertices are present is illustrated schematically in~\cref{fig:candidate-region}.

\begin{proposition}\label{prop:distinguished-degree}
    The distinguished vertices have degree at least $3$ in $R$.
\end{proposition}

\begin{proof}
    If $c_0$ is present, then it has at least $2$ neighbors in the $2$-connected block $R_B$ and at least $1$ neighbor outside $R_B$ in $R$, so $d(c_0) \geq 3$.

    If $u_0$ and $v_0$ are present, then $u_0$ is a cutvertex of $R_B-v_0$, since $R_B-v_0$ is connected but $R_B-u_0-v_0$ is disconnected.
    Hence, $u_0$ has at least $2$ neighbors in $R_B-v_0$ and is also adjacent to $v_0$, so $d(u_0)\geq3$.
    By symmetry, $d(v_0)\geq3$.
\end{proof}

\begin{proposition}\label{prop:distinguished-componentconnection}
    If $u_0$ and $v_0$ are present, then each component of $R-u_0-v_0$ is joined to both $u_0$ and $v_0$.
\end{proposition}

\begin{proof}
    By the $2$-connectivity of $R_B$, each component of $R_B-u_0-v_0$ is joined to both $u_0$ and $v_0$.
    The components of $R_B-u_0-v_0$ and $R-u_0-v_0$ are the same, except the one containing $c_0$ if $c_0$ is present, since in that case $c_0$ is distinct from $u_0$ and $v_0$.
    If $c_0$ is present, the component containing $c_0$ in $R-u_0-v_0$ is exactly the component containing $c_0$ in $R_B-u_0-v_0$ attached to $R-(V(R_B) \setminus \{c_0\})$ via $c_0$.
    Hence, each component of $R-u_0-v_0$ is connected to both $u_0$ and $v_0$, regardless of whether $c_0$ exists.
\end{proof}

We highlight the reason behind the choice of the distinguished vertices and $D$ with the following lemma.
\begin{lemma}\label{lem:region}
    Let $uv$ be a problematic edge-bundle in $D$ such that $\{u,v\} \neq \{u_0,v_0\}$ (if $u_0v_0$ exists) and $u,v \neq c_0$ (if $c_0$ exists).
    Then $R - u - v$ is connected.
\end{lemma}

\begin{proof}
    First, assume that $u_0v_0$ does not exist.
    Then, by construction, $uv$ cannot be an $R_B$-separating problematic edge-bundle; that is, $R_B-u-v$ is connected.
    If $c_0$ does not exist, the statement immediately follows, as $R = R_B$.
    Otherwise, since $c_0 \notin \{u,v\}$, the part of $R$ outside $R_B-c_0$ remains untouched and is attached to the connected graph $R_B-u-v$ via $c_0$.
    Hence, $R-u-v$ is connected.

    Now assume that $u_0v_0$ exists.
    By symmetry, we may assume that $u_0 \notin \{u,v\}$.
    Since $uv$ is an edge-bundle in $D$, we have $\{u,v\} \subseteq V(Q_0) \cup \{v_0\}$, and at least one of $u,v$ belongs to $V(Q_0)$.
    Every component of $R-u_0-v_0$ distinct from $Q_0$ remains intact and is joined to $u_0$ by~\cref{prop:distinguished-componentconnection}.
    If $v_0 \notin \{u,v\}$, then it also remains reachable from $u_0$ via the edge $u_0v_0$.
    Consequently, if $R-u-v$ were disconnected, then a component not containing $u_0$ would have its vertex set strictly contained in $V(Q_0)$.
    This would also be a component of $R_B-u-v$ not containing $c_0$ (if $c_0$ exists), contradicting the minimality of $|V(Q_0)|$ in the choice of $u_0v_0$ and $Q_0$.
    Therefore, $R-u-v$ is connected.
\end{proof}

Henceforth, write $A(w)=A_G(w)$ for every vertex $w \in V(R_B)$.
For every $w\in V(R_B)$ distinct from $c_0$ if $c_0$ is present, every edge of $G$ incident with $w$ belongs to $B$.
If a neighbor $v$ of $w$ in $G$ is deleted during compression, then $w$ is an endpoint of a suppressed path, so $d_G(v)=2<d_G(w)$ and $v\notin A(w)$.
For such a vertex $w$, it follows from the definition of problematic bundles in $R$ and the preservation of degrees that $A(w)$ consists precisely of the neighbors $v$ of $w$ in $R_B$ for which $wv$ is problematic in $R$ and, in the opposite-sign case, $d(v) \geq d(w)$.
Define a set of vertices $W$ in $D$ as
\[
W \coloneqq
\begin{cases}
    V(Q_0)&\text{if } u_0v_0 \text{ exists,}\\
    \{w \in V(R_B): d(w) \geq 3, w \neq c_0, c_0 \notin A(w)\}&\text{if only } c_0 \text{ exists,}\\
    V(R_B)&\text{otherwise.}
\end{cases}
\]

\begin{lemma}\label{lem:candidates}
    The set $W$ satisfies the following properties.
    \begin{enumerate}[label=(\roman*)]
        \item\label{cand:i} $W$ is non-empty.
        \item\label{cand:ii} Each vertex of $W$ has degree at least $3$.
        \item\label{cand:iii} Let $w \in W$.
          Then $w$ is distinct from all distinguished vertices.
          Furthermore, $c_0$ is not in $A(w)$ (if it exists), and at most one of $u_0$ and $v_0$ is contained in $A(w)$ (if they exist).
        \item\label{cand:iv} If $w \in W$ and $v \in A(w)$, then $R - w - v$ is connected.
        In other words, $v$ is not a cutvertex in $R-w$.
        \item\label{cand:v} The distinguished vertices (provided that at least one of them exists) not contained in $A(w)$ lie in the same component of $R-w - A(w)$ for each $w \in W$.

    \end{enumerate}
\end{lemma}

\begin{proof}
    We prove the properties in this order by case analysis, following the cases in the definition of $W$.

    \smallskip
    \noindent\emph{Proof of~\cref{cand:i}.}
    If $u_0v_0$ exists, $W = V(Q_0) \neq \emptyset$ is obvious, as $Q_0$ is a component of $R_B-u_0-v_0$.
    If no distinguished vertex exists, then $W=V(R_B) \neq \emptyset$.
    It remains when only $c_0$ exists.

    Suppose, for a contradiction, that $W = \emptyset$, that is, for each vertex $w \neq c_0$ of degree at least $3$ in $R_B$, the vertex $c_0$ is contained in $A(w)$.
    If there is a suspended path of length $2$ or $4$ incident with $c_0$ in $R_B$, then the other endpoint $w$ of that path has degree at least $3$ and is in the same bipartition class as $c_0$.
    In this case, we have $w \in W$, since $w$ cannot be adjacent to $c_0$, reaching a contradiction.
    Therefore, there is no suspended path of length $2$ or $4$ incident with $c_0$ in $R_B$, so each vertex in $R_B$ has degree at least $3$.
    Consequently, every vertex of $R_B$ distinct from $c_0$ is adjacent to $c_0$.
    Since $R_B$ contains at least $3$ vertices, there are at least $2$ neighbors; however, this contradicts the $2$-connectivity of $R_B$, since these neighbors are in the same bipartition class.

    \smallskip
    \noindent\emph{Proof of~\cref{cand:ii}.}
    If $c_0$ is absent, then all vertices of $R_B$ have $R$-degree at least $3$, as observed after compression.
    If $c_0$ is the only distinguished vertex, then the assertion follows from the definition of $W$.
    It remains to consider the case where all three distinguished vertices are present.
    The vertices $u_0$ and $v_0$ are distinct from $c_0$ by definition and have degree at least $3$ by~\cref{prop:distinguished-degree}.
    Thus, every vertex of $R_B$ with $R$-degree $2$ lies in the component of $R_B-u_0-v_0$ containing $c_0$, since it is an internal vertex of a suspended path of length $2$ or $4$ incident with $c_0$.
    As this component is distinct from $Q_0$, every vertex of $W=V(Q_0)$ has degree at least $3$.

    \smallskip
    \noindent\emph{Proof of~\cref{cand:iii}.}
    First, assume that the vertices $u_0$ and $v_0$ are present.
    In this case, we have $W = V(Q_0)$.
    None of the distinguished vertices is contained in $Q_0$, as $Q_0$ is a component of $R-u_0-v_0$ not containing $c_0$ if $c_0$ is present.
    Hence, $w$ is distinct from all distinguished vertices.
    Furthermore, we have $A(w) \subseteq N_R(w) \subseteq V(Q_0) \cup \{u_0,v_0\}$.
    If $c_0$ is present, it follows that $c_0 \notin A(w)$, as $c_0$ is distinct from $u_0$ and $v_0$ by definition.
    Since $A(w)$ is an independent set, it may contain at most one of the adjacent pair $u_0,v_0$.

    If $c_0$ is the only distinguished vertex that exists, the statement follows immediately from the definition of $W$.
    If no distinguished vertex is present, the statement holds vacuously.

    \smallskip
    \noindent\emph{Proof of~\cref{cand:iv}.}
    It follows immediately from~\cref{cand:iii} that $\{w,v\} \neq \{u_0,v_0\}$ if $u_0,v_0$ are present, and neither $w$ nor $v$ is $c_0$ if $c_0$ is present.
    Hence, the desired statement follows from~\cref{lem:region}, as $wv$ is problematic by definition.

    \smallskip
    \noindent\emph{Proof of~\cref{cand:v}.}
    If no distinguished vertex is present or if $c_0$ is the only distinguished vertex, the statement holds immediately.
    Next, assume that $c_0$ does not exist, although $u_0$ and $v_0$ are present.
    If $u_0$ or $v_0$ is contained in $A(w)$, then the statement is obvious.
    Otherwise, $u_0$ and $v_0$ must be in the same component, since they are adjacent in $R$.

    The remaining case is where each possible distinguished vertex is present.
    Then we have $W= V(Q_0)$ and $A(w) \subseteq V(Q_0) \cup \{u_0,v_0\}$.
    By~\cref{cand:iii}, we have $c_0 \notin A(w)$ and $u_0$ or $v_0$ is not contained in $A(w)$.
    By symmetry, we may assume that $u_0 \notin A(w)$ and deduce that $A(w) \subseteq V(Q_0) \cup \{v_0\}$.
    By~\cref{prop:distinguished-componentconnection}, each component of $R-u_0-v_0$ is connected to both $u_0$ and $v_0$.
    Therefore, the vertices outside $V(Q_0) \cup \{v_0\}$ are in the same component of $R-w-A(w)$, and, if $v_0 \notin A(w)$, then $v_0$ is also in this component.
\end{proof}

\subsection{Selecting a candidate vertex}

We define a subset $W^*\subseteq W$ as follows.
If $a,b$ have the same sign, then let $W^*$ be the set of vertices $w \in W$ that are local extrema among their neighbors of degree at least $3$ in $G$.
If $a,b$ have opposite signs, then put $W^* = W$.

\begin{proposition}
    $W^*$ is non-empty.
\end{proposition}

\begin{proof}
    By~\cref{lem:candidates}\cref{cand:i}, the set $W$ is non-empty, so the assertion holds in the opposite-sign case, where $W^*=W$.
    Henceforth, assume that $a,b$ have the same sign.
    Every neighbor $v$ in $G$ of a vertex $w \in W$ with $d_G(v)\geq3$ is also a neighbor of $w$ in $R_B$.
    Since the degrees of retained vertices are unchanged, every vertex of $W$ that is a local extremum among its neighbors in $R$ belongs to $W^*$.

    If no distinguished vertex exists, then $W=V(R_B)=V(R)$, so a vertex of maximum degree in $W$ belongs to $W^*$.

    In the remaining cases, define
    \[
        U \coloneqq
        \begin{cases}
            \{u_0,v_0\}& \text{if } u_0,v_0 \text{ are present,}\\
            \{c_0\} \cup A(c_0)& \text{if only } c_0 \text{ is present.}
        \end{cases}
    \]
    In the first case, $W=V(Q_0)$, so every neighbor in $R$ of a vertex of $W$ belongs to $W \cup U$.
    In the second case, every vertex of $R_B$ has $R$-degree at least $3$, and the definition of problematic bundles gives $c_0 \in A(w)$ if and only if $w \in A(c_0)$ for every $w \in V(R_B)\setminus\{c_0\}$.
    Hence, $W=V(R_B)\setminus U$, so again every neighbor in $R$ of a vertex of $W$ belongs to $W \cup U$.

    In both cases, the definition of a problematic bundle implies that all vertices of $U$ have the same degree, say $d_U$.
    If $W$ contains a vertex of degree at least $d_U$, then a vertex of maximum degree in $W$ is a local maximum among its neighbors in $R$.
    Otherwise, a vertex of minimum degree in $W$ is a local minimum among its neighbors in $R$.
    In either case, the selected vertex belongs to $W^*$.
\end{proof}

\paragraph{The choice rule.} We choose the vertex $w_0$ from $W^*$ according to the first applicable case below.
\begin{enumerate}
    \item\label{choice:case1} If there exists a vertex $w \in W^*$ with $|A(w)| \leq 1$, then choose such a vertex to be $w_0$.
    \item\label{choice:case2} If there exists a vertex of $W^*$ of degree at least $5$, then choose such a vertex to be $w_0$.
    \item\label{choice:case3} If neither of the previous two possibilities holds, let $C(w)$ denote the component of $R-w-A(w)$ that contains the distinguished vertices not contained in $A(w)$ (by~\cref{lem:candidates}\cref{cand:iii,cand:v}, this component is well-defined), or, if no distinguished vertex is present, an arbitrarily chosen component of maximum order.
    If no distinguished vertex is present, then $R=R_B$ is $2$-connected with at least three vertices, so the bipartition class of $w$ contains another vertex, which survives the deletion of $w$ and $A(w)$.
    Choose $w_0 \in W^*$ such that $|V(C(w_0))|$ is maximum over all such choices.
\end{enumerate}

When Case~\cref{choice:case3} applies, write $C_0=C(w_0)$.
In this case, for every $w \in W^*$, no component of $R-w-A(w)$ has a vertex set strictly containing $V(C_0)$.
If such a component existed and distinguished vertices were present, it would be $C(w)$, since $C_0$ contains at least one distinguished vertex by~\cref{lem:candidates}\cref{cand:iii,cand:v}.
If no distinguished vertex were present, then $C(w)$ would have at least as many vertices as this component.
In either case, we would have $|V(C(w))|>|V(C_0)|$, contradicting the choice of $w_0$.

Next, we establish some properties of $w_0$ in Case~\cref{choice:case3}.

\begin{lemma}\label{lem:ext-struct}
    Assume that Case~\cref{choice:case3} applies above.
    Then the following conditions hold, where in~\cref{ext-struct:ii,ext-struct:iii,ext-struct:iv,ext-struct:v}, $C$ denotes an arbitrary component of $R-w_0-A(w_0)$ distinct from $C_0$.
    \begin{enumerate}[label=(\roman*)]
        \item\label{ext-struct:i}
        If $c_0$ is present, then $c_0$ and the internal vertices of suspended paths of length $2$ or $4$ incident with $c_0$ are contained in $V(C_0)$.
        In that case, each vertex connected to $c_0$ by an edge or a suspended path in $R$ is in $V(C_0)$ or is contained in $A(w_0) \cup \{w_0\}$ and has an edge to $C_0$.
        If $u_0,v_0$ are present, then both are in $V(C_0) \cup A(w_0)$ with at most one of them in $A(w_0)$, and, if one of them is in $A(w_0)$, then it has an edge to $C_0$.

        \item\label{ext-struct:ii} Each vertex of $C$ has degree at least $3$.
        \item\label{ext-struct:iii} Each vertex of $C$ is contained in $W$.
        \item\label{ext-struct:iv} $C$ only contains vertices of the same degree as $w_0$.
        \item\label{ext-struct:v} $C$ contains exactly one vertex, and that vertex has the same degree as $w_0$.
        \item\label{ext-struct:vi} There exists an edge from each vertex of $A(w_0)$ to $C_0$.
    \end{enumerate}
\end{lemma}

\begin{proof}
    Since Cases~\cref{choice:case1} and~\cref{choice:case2} do not apply, we have $|A(w_0)|\geq2$ and $d(w_0)\in\{3,4\}$.
    Since $R-w_0$ is connected and $A(w_0)$ is non-empty, some vertex of $A(w_0)$ has an edge to $C_0$.

    \smallskip
    \noindent\emph{Proof of~\cref{ext-struct:i}.}
    Suppose first that $c_0$ is present.
    Then $c_0 \in V(C_0)$ by the choice of $C_0$ and the fact $c_0 \notin A(w_0)$ provided by~\cref{lem:candidates}\cref{cand:iii}.
    Since $d(w_0)\geq3$ and $d(v)\geq d(w_0)$ for every $v\in A(w_0)$, none of the deleted vertices is an internal vertex of a suspended path incident with $c_0$.
    Thus, the internal vertices of the suspended paths incident with $c_0$ in $R$ remain in the component containing $c_0$, which is $C_0$.
    Furthermore, if a vertex $v \notin A(w_0) \cup \{w_0\}$ is connected to $c_0$ by an edge or a suspended path in $R$, then that edge or path remains in $R-w_0-A(w_0)$.
    This implies that $c_0$ and $v$ must be in the same component, which is $C_0$.
    If $v \in A(w_0) \cup\{w_0\}$ is connected to $c_0$ by an edge or a suspended path in $R$, then there exists an edge from $v$ to $c_0$ or to an internal vertex of a suspended path incident with $c_0$; both are contained in $V(C_0)$.

    If $u_0$ and $v_0$ are present, then by~\cref{lem:candidates}\cref{cand:iii,cand:v}, both belong to $V(C_0)\cup A(w_0)$, and at most one belongs to $A(w_0)$.
    If one of them belongs to $A(w_0)$, the bundle $u_0v_0$ connects that vertex to $C_0$.

    \smallskip
    \noindent\emph{Proof of~\cref{ext-struct:ii}.}
    If $c_0$ is present, then all deleted vertices belong to $V(R_B)\setminus\{c_0\}$, so all vertices of $R$ outside $R_B$ remain connected to $c_0$ and hence belong to $C_0$.
    Therefore, $C$ is contained in $R_B$.
    By the $2$-connectivity of $R_B$, each vertex of $C$ has $R$-degree at least $2$.
    Vertices of $R_B$ with $R$-degree $2$ may only occur as internal vertices of suspended paths of length $2$ or $4$ incident with $c_0$.
    These are contained in $C_0$ by~\cref{ext-struct:i}, which proves~\cref{ext-struct:ii}.

    \smallskip
    \noindent\emph{Proof of~\cref{ext-struct:iii}.}
    First, assume that $u_0$ and $v_0$ are present.
    Then at most one of $u_0$ and $v_0$ is contained in $A(w_0)$ by~\cref{lem:candidates}\cref{cand:iii}.
    By symmetry, we may assume that $u_0 \notin A(w_0)$.
    We have $\{w_0\} \cup A(w_0) \subseteq V(Q_0) \cup \{v_0\}$, so each surviving vertex outside $V(Q_0)$ must be in the same component as $u_0$ in $R-w_0-A(w_0)$ by~\cref{prop:distinguished-componentconnection}, which is exactly $C_0$.
    Hence, $V(C) \subseteq V(Q_0) = W$.

    Now assume that only $c_0$ is present.
    Then the vertices $v\in V(R_B)$ of degree at least $3$ with $c_0 \in A(v)$ are either contained in $A(w_0) \cup \{w_0\}$ or in $C_0$ by~\cref{ext-struct:i}.
    The internal vertices of suspended paths incident with $c_0$ are also in $C_0$.
    Therefore, we have $V(C) \subseteq W$ by the definition of $W$.

    If no distinguished vertex exists, we have $W = V(R_B)$ by definition, so the statement is obvious.

    \smallskip
    \noindent\emph{Proof of~\cref{ext-struct:iv}.}
    Assume to the contrary that $C$ contains at least one vertex of different degree from that of $w_0$.

    \smallskip
    \noindent\textbf{Possibility 1.} The weights $a,b$ have the same sign.\\
    In this case, each vertex of $A(w_0)$ has degree $d(w_0)$.
    Assume that there exists a vertex in $C$ of higher degree than that of $w_0$ and choose such a vertex $w$ of maximum degree in $C$.
    Since $C$ is a component of $R-w_0 -A(w_0)$, we have $N_R(w) \subseteq V(C) \cup \{w_0\} \cup A(w_0)$.
    Hence, the vertex $w$ is a local maximum in $R$, so $w \in W^*$.
    Furthermore, $A(w) \subseteq V(C)$ holds, as $w$ has different degree from $d(w_0)$, which is the common degree of the vertices in $\{w_0\} \cup A(w_0)$.
    It follows that $R-w-A(w)$ has a component that contains the vertices of $C_0$ and at least one vertex in $\{w_0\} \cup A(w_0)$, contrary to the choice of $w_0$.
    If $C$ contains a vertex of lower degree than that of $w_0$, we arrive at the same conclusion if we choose $w$ as the minimum degree vertex in $C$.
    This proves the statement when the weights have the same sign.

    \smallskip
    \noindent\textbf{Possibility 2.} The weights $a,b$ have opposite signs.\\
    In this case, we have $V(C) \subseteq W = W^*$.
    Therefore, $C$ only contains vertices of degree $3$ or $4$, as Case~\cref{choice:case2} does not apply.
    Note that the degrees in $A(w_0)$ are congruent to $d(w_0)$ modulo~$\delta$.
    Assume, for a contradiction, that $C$ contains a vertex $w$ of different degree from that of $w_0$.
    Since $w$ and $w_0$ both have degree $3$ or $4$, this implies that $d(w) \not\equiv d(w_0) \pmod{\delta}$.
    We can argue as in the same-sign case.
    We have $N_R(w) \subseteq V(C) \cup \{w_0\} \cup A(w_0)$, and therefore $A(w) \subseteq V(C)$, as $d(w)$ is incongruent to each degree in $\{w_0\} \cup A(w_0)$ modulo~$\delta$.
    This also contradicts the choice of $w_0$ by the same argument as in the same-sign case, proving the statement in the opposite-sign case.

    \smallskip
    \noindent\emph{Proof of~\cref{ext-struct:v}.}
    Assume that $C$ contains more than one vertex.
    Since $R$ is bipartite and $C$ is connected, the component $C$ must contain a vertex $w$ in the same bipartition class as the vertices of $A(w_0)$.
    This implies that $A(w) \subseteq N_R(w) \subseteq V(C) \cup \{w_0\}$.
    By~\cref{ext-struct:iv}, each vertex in $V(C) \cup \{w_0\}$ has the same degree as $w_0$.
    It follows that $w$ is a local extremum, so $w \in W^*$ holds regardless of the signs of the weights.
    This contradicts the choice of $w_0$, as $R-w-A(w)$ has a component containing the vertices of $C_0$ and at least one vertex of $A(w_0)$.

    \smallskip
    \noindent\emph{Proof of~\cref{ext-struct:vi}.}
    Suppose, for a contradiction, that there exists a vertex $w \in A(w_0)$ without an edge to $C_0$.
    Then, by~\cref{ext-struct:i}, the vertex $w$ is not a distinguished vertex and, if $c_0$ is present, is not connected to $c_0$ by an edge or a suspended path.

    We claim that $w \in W^*$.
    We first show that $w \in W$.
    If $u_0,v_0$ are present, then $w \in V(Q_0)=W$, since $w$ is a neighbor of $w_0$ distinct from $u_0,v_0$.
    If only $c_0$ is present, then $w$ satisfies the requirements of $W$, since $d(w) \geq 3$, $w \neq c_0$, and $c_0 \notin N_R(w)$.
    If no distinguished vertex is present, the claim is obvious because $W = V(R_B)$.
    By our assumption, $w$ has no neighbors in $V(C_0)$; furthermore, $w$ has no neighbors in $A(w_0)$ either, as the vertices of $A(w_0)$ are in the same bipartition class.
    Hence, the neighbors of $w$ in $R$ distinct from $w_0$ lie in the components distinct from $C_0$.
    In the same-sign case, $d(w)=d(w_0)$ by the definition of $A(w_0)$, and by~\cref{ext-struct:iv}, the vertices of these components have the same degree as $w_0$, so $w$ is a local extremum.
    Thus, $w \in W^*$ in both sign cases.
    Since $A(w)\subseteq N_R(w)$, neither a vertex of $C_0$ nor a vertex of $A(w_0)$ other than $w$ is deleted when forming $R-w-A(w)$.
    As some vertex of $A(w_0)$ has an edge to $C_0$ and $w$ has no such edge, a component of $R-w-A(w)$ contains all vertices of $C_0$ and at least one vertex of $A(w_0)$.
    This contradicts the choice of $w_0$ and proves~\cref{ext-struct:vi}.
\end{proof}

\subsection{Validity of the choice}

We shall now verify that $w_0$ is indeed an enforceable extremum, except if a specific cubic configuration occurs.
The following hypotheses of~\cref{lem:samesign-locminmax,lem:opsign-locmax} are already satisfied.
\begin{itemize}
    \item The vertex $w_0$ is not a cutvertex of $G$, as it belongs to $V(B)$ and is distinct from $c_0$ if $c_0$ is present.
    \item The degree of $w_0$ is at least $3$ by~\cref{lem:candidates}\cref{cand:ii}.
    \item In the same-sign case, $w_0$ is a local extremum among its neighbors of degree at least $3$ in $G$ according to the definition of $W^*$.
    \item In the same-sign case, each suspended path incident with $w_0$ in $G$ belongs to $B$ and has length $3$, as established before compression.
\end{itemize}
Therefore, it remains to show that there exists a set of edges $F \subseteq E(G^*(w_0))$ that satisfies the corresponding lower bounds on the degrees in $A(w_0)$ and leaves a connected graph upon deletion.

Notice that it suffices to show that such a set exists in $E(R-w_0)$.
To see this, assume that we found an appropriate set of edges $F_R$ in $R-w_0$.
We may assume that $F_R$ only contains edges with an endpoint in $A(w_0)$, since the removal of other edges from $F_R$ does not change the $F_R$-degrees in $A(w_0)$, and $R-w_0-F_R$ remains connected, since edges are only restored.
Since $A(w_0)$ is independent, each blue edge in $F_R$ has exactly one endpoint in $A(w_0)$.
Replace each such blue edge by the outer edge of the corresponding suspended path of length $3$ incident with that endpoint to obtain a set of edges $F \subseteq E(G^*(w_0))$.
The graph $G^*(w_0)$ is obtained from $R-w_0$ by expanding its blue edges into the corresponding suspended paths.
Consequently, $G^*(w_0)-F$ is obtained from the connected graph $R-w_0-F_R$ by subdividing edges and attaching paths of length $2$, so it is connected.
The $F$-degrees and $F_R$-degrees in $A(w_0)$ are the same, so $F$ is an appropriate set of edges.

We will show that $F_R$ exists using~\cref{cor:frank-leaving}.
Abbreviate $A = A(w_0)$, write $d = d(w_0)$ and $J = R-w_0$.
Since $w_0$ is not a cutvertex in $R$, the graph $J$ is connected.
Introduce the notation $p = |a| + |b| = a+ |b|$.
Since $|b|\geq a+1$, we have $p\geq 2a+1\geq 3$.
Also note that in the opposite-sign case, $p = a-b = \delta$.
Define $k_x= \frac{1}{p}(d(x) - d(w_0))$ for each $x \in A(w_0)$.
Define $h\colon A \to \Z_{\geq 0}$ as
\[
h(x) \coloneqq \frac{|a|}{p}(d(x) - d(w_0)) = a k_x.
\]
With these notations, we can handle the same-sign and opposite-sign cases together.
By the definition of $A(w_0)$, each $x \in A$ satisfies
\begin{itemize}
    \item $k_x = h(x) = 0$ in the same-sign case, as $d(x) = d(w_0)$,
    \item $k_x, h(x) \in \Z_{\geq 0}$ in the opposite-sign case, since $d(x) \geq d(w_0)$ and $d(x) \equiv d(w_0) \pmod{p}$ hold, as $p = \delta$.
\end{itemize}
In~\cref{lem:samesign-locminmax,lem:opsign-locmax}, the lower bound on the degree of $x$ is exactly $h(x) + 1$ in both cases: in the same-sign case, $h(x) = 0$; in the opposite-sign case, $h(x) = \frac{|a|}{\delta}(d(x) - d(w_0))$.
These are non-negative integer bounds, the graph $J$ is connected, and the set $A$ is independent, so~\cref{cor:frank-leaving} applies.
Therefore, it suffices to verify that
\begin{equation}\label{eq:frank-leaving-special}
\sum_{x \in X}d_J(x) - a \sum_{x \in X} k_x \geq 2|X|+\sigma_J(X)-1
\end{equation}
holds for each subset $\emptyset \neq X \subseteq A$, since $d_J(x) - h(x) = d_J(x) -ak_x$ by definition.

Let $m_x$ denote the multiplicity of the bundle $w_0x$ in $R$ for each $x \in A$.
Note that we have $m_x \geq 1$ and $\sum_{x \in A} m_x \leq d$.
We can express $d_J(x)$ in terms of $d, k_x$, and $m_x$ as
\begin{equation}\label{eq:degree-expression}
    d_J(x) = d_{R-w_0}(x) = d(x) - m_x = d + (d(x) - d) - m_x = d + p k_x - m_x.
\end{equation}
Since $X$ is independent, we have $d_J(X) = \sum_{x \in X} d_J(x)$.
Write
\[
K_X = \sum_{x \in X} k_x, \qquad M_X = \sum_{x \in X} m_x.
\]
Using this notation and~\cref{eq:degree-expression}, inequality~\cref{eq:frank-leaving-special} can be rewritten equivalently as
\begin{equation}\label{eq:frank-leaving-special2}
    d |X| + (p-a)K_X - M_X \geq 2 |X| + \sigma_J(X) -1.
\end{equation}

\begin{lemma}\label{lem:modified-deg-bound}
    Each vertex $x \in A$ satisfies $d + \ell k_x -m_x \geq \max \{2, 1+k_x\}$ for any positive integer $\ell$.
\end{lemma}

\begin{proof}
    If $m_x \leq d-2$, the statement follows immediately, as $\ell k_x \geq k_x \geq 0$.
    Suppose, for a contradiction, that $m_x = d$.
    Then $w_0$ has only one neighbor $x$ in $R$; therefore, the deletion of $x$ would isolate $w_0$, while vertices in $R_B$ distinct from $w_0$ and $x$ remain, since $R_B$ has at least $3$ vertices.
    This contradicts the $2$-connectivity of $R_B$.

    The remaining case is $m_x = d-1$.
    If $k_x \geq 1$, then $\ell k_x \geq k_x \geq 1$, and the statement follows.
    If $k_x = 0$, then $d(x) = d(w_0)$, so the bundle $w_0 x$ only contains white edges according to the definition of $A(w_0)$.
    This results in exactly the configuration excluded by~\cref{prop:minimality-prop}\cref{min:M2}: two vertices of degree $q+1 = d$ connected by $q = d-1 \geq 2$ edges in the original graph $G$.
    We covered all possible cases; the proof is complete.
\end{proof}

\begin{corollary}\label{cor:deg-geq2}
    Each vertex $x \in A$ satisfies $d_J(x) \geq 2$.
\end{corollary}

\begin{proof}
    By~\cref{eq:degree-expression}, we have $d_J(x) = d + p k_x - m_x$, so the statement follows from~\cref{lem:modified-deg-bound} with $\ell = p$.
\end{proof}

\begin{lemma}\label{lem:frank-small}
    Inequality~\cref{eq:frank-leaving-special} holds for each $X \subseteq A$ with $|X| =1$.
\end{lemma}

\begin{proof}

    Assume that $X = \{x\}$.
    Since~\cref{eq:frank-leaving-special} and~\cref{eq:frank-leaving-special2} are equivalent, it suffices to verify the latter.
    Since $x \in A(w_0)$, the graph $J-x = R-w_0-x$ is connected by~\cref{lem:candidates}\cref{cand:iv}.
    Therefore, we have $\sigma_J(X) =1$, and the value of the right-hand side is $2$.
    The left-hand side takes the form $d+(p-a)k_x -m_x$.
    Hence, the desired inequality follows from~\cref{lem:modified-deg-bound} with $\ell = p-a = |b| > 0$.
\end{proof}

\begin{corollary}\label{cor:enfext-case1}
    The vertex $w_0$ is an enforceable extremum if Case~\cref{choice:case1} applies.
\end{corollary}

\begin{proof}
    If Case~\cref{choice:case1} applies, then $|A| \leq 1$.
    Hence, each non-empty subset $X \subseteq A$ satisfies the desired inequality by~\cref{lem:frank-small}.
\end{proof}

\begin{lemma}\label{lem:deg-compnum-ineq}
    Each subset $\emptyset \neq X \subseteq A$ satisfies $d_J(X) \geq 2\sigma_J(X)$.
\end{lemma}

\begin{proof}
    Fix an arbitrary non-empty subset $X$.
    We use two properties of the elements of $X$.
    By~\cref{lem:candidates}\cref{cand:iv}, no vertex of $X$ is a cutvertex in $J=R-w_0$, and by~\cref{cor:deg-geq2}, each vertex of $X$ has $J$-degree at least $2$.

    Let $K$ be a component of $J-X$.
    We claim that there are at least $2$ edges between $K$ and $X$ in $J$.
    Suppose, for a contradiction, that $K$ and $X$ are connected by at most $1$ edge.
    Since $J$ is connected, they are connected by precisely one edge; let $x$ be the unique vertex of $X$ with an edge to $K$.
    We have $|X| \geq 2$, since otherwise $d_J(x) = 1$ would hold, a contradiction.
    Hence, $x$ is a cutvertex in $J$, since it separates $K$ and $X \setminus \{x\}$.
    This is also a contradiction, which proves our claim.

    Since $d_J(X)$ counts the edges between $X$ and $V(J) \setminus X$ and each component of $J-X$ is connected to $X$ by at least two edges, we obtain $d_J(X) \geq 2 \sigma_J(X)$, which is the desired inequality.
\end{proof}

\begin{lemma}\label{lem:mulsum-bound}
    Each subset $ \emptyset \neq X \subseteq A$ satisfies $M_X \leq d+ |X| -|N_R(w_0)|$.
\end{lemma}

\begin{proof}
    Let $r$ denote the number of edges in $R$ incident with $w_0$ with an endpoint not in $X$.
    Then, we have $r \geq |N_R(w_0) \setminus X| = |N_R(w_0)| - |X|$, as each neighbor of $w_0$ other than the vertices of $X$ has an edge to $w_0$.
    Therefore, since $M_X$ counts the edges from $w_0$ to $X$, we have
    \[
    M_X = d - r  \leq d-(|N_R(w_0)| - |X|) = d+ |X| - |N_R(w_0)|,
    \]
    which is the desired inequality.
\end{proof}

Introduce the notation $\eta(X) = d_J(X) - 2 \sigma_J(X)$.
By~\cref{lem:deg-compnum-ineq}, we have $\eta(X) \geq 0$ for any subset $\emptyset \neq X \subseteq A$.

\begin{proposition}\label{prop:enfext-case2}
    The vertex $w_0$ is an enforceable extremum if Case~\cref{choice:case2} applies.
\end{proposition}

\begin{proof}
    In Case~\cref{choice:case2}, we have $d=d(w_0)\geq5$.
    By~\cref{lem:frank-small}, it suffices to verify inequality~\cref{eq:frank-leaving-special} for subsets $X\subseteq A$ with $|X|\geq2$.
    Fix such a subset and let $S=d_J(X)-aK_X-2|X|-\sigma_J(X)+1$ be the integer slack in that inequality.
    Using the definition of $\eta(X)$,~\cref{eq:degree-expression}, and the bound in~\cref{lem:mulsum-bound}, we obtain
    \[
    2S = d|X|+(p-2a)K_X-M_X+\eta(X)-4|X|+2 \geq (d-5)(|X|-1)+(p-2a)K_X+\eta(X)+|N_R(w_0)|-3.
    \]
    Since $d\geq5$, $p>2a$, and $K_X,\eta(X)\geq0$, the first three terms of the lower bound are non-negative.
    Moreover, $|N_R(w_0)|\geq|X|\geq2$, so $2S\geq-1$.
    As $S$ is an integer, we have $S\geq0$, as desired.
\end{proof}

\begin{lemma}\label{lem:components-case3}
    Assume that Case~\cref{choice:case3} applies.
    Then for each $\emptyset \neq X \subseteq A$, the graph $J-X$ has a component connected to each vertex of $X$ by an edge, and all other components are isolated vertices of degree $d$ in $J$.
\end{lemma}

\begin{proof}
    By~\cref{lem:ext-struct}\cref{ext-struct:v,ext-struct:vi}, the components of $J-A$ are $C_0$ and singletons of degree $d$ in $R$, and every vertex of $A$ has an edge to $C_0$.
    Since $J$ is connected, each such singleton has a neighbor in $A$ and hence belongs to the same bipartition class as $w_0$ in $R$.
    Therefore, its $J$-degree equals its $R$-degree.

    For an arbitrary non-empty subset $X\subseteq A$, obtain $J-X$ from $J-A$ by restoring all vertices of $A\setminus X$ and their edges in $J-X$.
    Each restored vertex has an edge to $C_0$, so the component containing $C_0$ also contains every restored vertex and exactly those singleton components of $J-A$ with a restored neighbor.
    This component is connected to every vertex of $X$ by an edge, and the remaining components are unchanged singletons of degree $d$ in $J$.
\end{proof}

\begin{corollary}\label{cor:deg-compnum-ineq-case3}
    If Case~\cref{choice:case3} applies, then each $\emptyset \neq X \subseteq A$ satisfies $d_J(X) \geq d(\sigma_J(X)-1) + |X|$.
\end{corollary}

\begin{proof}
    The designated component of $J-X$ in the sense of~\cref{lem:components-case3} is connected to each vertex of $X$ by at least one edge, so there are at least $|X|$ edges between $X$ and this component.
    The other components of $J-X$ are singletons of degree $d$ in $J$.
    Since all edges incident with these singletons have their other endpoint in $X$, there are exactly $d$ edges between $X$ and one of these singleton components.
    Since the number of non-designated singleton components is $\sigma_J(X) -1$, we obtain the desired inequality.
\end{proof}

\begin{lemma}\label{lem:frank-small-case3}
    Assume that Case~\cref{choice:case3} applies.
    Then inequality~\cref{eq:frank-leaving-special} holds for each $\emptyset \neq X \subseteq A$ with $|X| \leq 2$.
\end{lemma}

\begin{proof}
    By~\cref{lem:frank-small}, it suffices to verify the statement for $|X| = 2$.
    Fix a subset $X = \{x,y\}$.
    By symmetry, we may assume that $d_J(x) \geq d_J(y)$.
    By~\cref{eq:degree-expression}, we have
    \begin{equation}\label{eq:distances}
    p(k_x-k_y) \geq m_x -m_y.
    \end{equation}
    Since $p \geq 3$ and $|m_x -m_y| \leq d-2 \leq 2$, this also implies $k_x \geq k_y$.
    Since neither $x$ nor $y$ is a cutvertex in $J$ by~\cref{lem:candidates}\cref{cand:iv}, each component of $J-X$ is connected to both $x$ and $y$ by at least one edge.
    Hence, we have $d_J(y) \geq \sigma_J(X)$.
    Introduce $\alpha = d_J(y) - \sigma_J(X) \geq 0$.
    Substituting $\sigma_J(X) = d_J(y) - \alpha$ and $|X| = 2$, and canceling out $d_J(y)$ on both sides of~\cref{eq:frank-leaving-special}, we get $d_J(x) - a(k_x + k_y) \geq 3 - \alpha$.
    By substituting $d_J(x) = d + pk_x -m_x$ and rearranging, we get
    \begin{equation}\label{eq:frank-leaving-two-element}
    (d + (p-2a) k_x - m_x) + a(k_x -k_y) + \alpha \geq 3.
    \end{equation}

    Notice that each term is non-negative; furthermore, the first term is always at least $2$ by~\cref{lem:modified-deg-bound} with $\ell = p-2a \geq 1$.
    If $k_x \geq 2$, then the first term is at least $3$ by the same lemma, and the statement follows.

    Therefore, the inequality may only fail if $\alpha = 0$, $k_x \leq 1$.
    Then we have $d_J(y) = \sigma_J(X)$, and hence we have $d_J(x) + d_J(y) \geq d(d_J(y) - 1) + |X|$ by~\cref{cor:deg-compnum-ineq-case3}.
    After rearranging, we obtain
    \[
    d_J(x) \geq (d-1)d_J(y) -d + 2 = (d-1)(d_J(y) -1) +1.
    \]
    If $k_x = 0$, then $d_J(x) \leq d-1$, and the inequality above implies $d_J(y) = 1$, which contradicts~\cref{cor:deg-geq2}.
    If $k_x = 1$, then the inequality may only fail if $k_y = k_x = 1$.
    In this case, we have $m_x \leq m_y$ by~\cref{eq:distances}, so $m_x \leq \left\lfloor\frac{d}{2}\right\rfloor \leq d-2$, as $d \geq 3$.
    Hence, the first term in~\cref{eq:frank-leaving-two-element} is at least $3$, and the inequality also holds.
    The proof is complete, as all possible cases are covered.
\end{proof}

\begin{proposition}\label{prop:enfext-case3}
    Assume that Case~\cref{choice:case3} applies.
    Then $w_0$ is an enforceable extremum, unless the following exceptional configuration occurs: $d = |A| = 3$, each vertex of $A$ has degree $2$ in $J$, and $J-A$ has two components, each of which is connected to each vertex of $A$ by exactly one edge.
\end{proposition}

\begin{proof}
    By~\cref{lem:frank-small-case3}, it suffices to verify inequality~\cref{eq:frank-leaving-special} for each $X \subseteq A$ with $|X| \geq 3$.
    Introduce the notation $\eta'(X) = d_J(X) - (d(\sigma_J(X) -1) + |X|)$, which is non-negative by~\cref{cor:deg-compnum-ineq-case3}.
    Multiplying~\cref{eq:frank-leaving-special} by $d$ and using the definition of $\eta'(X)$ and~\cref{eq:degree-expression}, we obtain the equivalent inequality
    \[
    (d^2-3d+1)|X|+((d-1)p-da)K_X+\eta'(X)\geq(d-1)M_X.
    \]
    Since $p\geq2a+1$ and $a\geq1$, we have $(d-1)p-da\geq(d-2)a+d-1\geq2d-3$.

    Since Case~\cref{choice:case3} applies, we have $d\in\{3,4\}$.
    If $d=4$, then the left-hand side is at least $5|X|\geq15$ and the right-hand side is at most $12$, since $M_X\leq d$.

    If $d=3$, then $|X|=3$, so $X=A=N_R(w_0)$, $M_X=3$, and $d_J(X)=6+pK_X$.
    The inequality reduces to $(2p-3a)K_X+\eta'(X)\geq3$ and therefore holds if $K_X\geq1$.
    Henceforth, assume that $K_X= 0$.
    Then $d_J(X) = 6$, so
    \[
    \eta'(X) = d_J(X) - |X| - 3(\sigma_J(X) - 1) = 6-3 \sigma_J(X) = 3(2-\sigma_J(X)) \geq 0.
    \]
    Hence, the desired inequality fails only if $\sigma_J(X) = 2$.
    Since $|X|=3$ and $d_J(X)=6$, each vertex of $X$ has $J$-degree $2$ by~\cref{cor:deg-geq2}.
    By~\cref{lem:components-case3}, one component of $J-X$ has at least one edge to each vertex of $X$ and the other has three edges to $X$, forcing exactly one edge from each vertex of $X$ to each component.
    This results in the exceptional case, as $X = A$.
    In this configuration, inequality~\cref{eq:frank-leaving-special} for $X=A$ reads $6\geq7$ and therefore fails.
\end{proof}

\subsection{The remaining cubic configuration}

To complete the proof of~\cref{thm:main}, it remains to handle the exceptional case in~\cref{prop:enfext-case3}.
That is, assume that $d = |A| = 3$, each vertex of $A$ has degree $2$ in $J$ and thus degree $3$ in $R$, and $J-A$ has two components, each of which is connected to each vertex of $A$ by exactly one edge.
Since $d=|A|=3$, each vertex of $A$ is joined to $w_0$ by exactly one edge, which is white by the definition of $A$.
By~\cref{lem:ext-struct}\cref{ext-struct:iii,ext-struct:v}, the component distinct from $C_0$ consists of a single vertex $w \in W$ of degree $3$.
Thus, each vertex of $A$ is joined to $w$ by exactly one edge.
Since $N_R(w)=A$ and each vertex of $A$ has degree $3$, we have $w \in W^*$.

Suppose, for a contradiction, that there is a blue edge incident with $w$, and let $x \in A$ denote its other endpoint.
Since the bundle $wx$ consists of a single blue edge, $w$ and $x$ are not adjacent in $G$, so $A(w) \subseteq A \setminus \{x\}$.
Since $x$ has an edge to $C_0$, the component of $R-w-A(w)$ containing $x$ also contains all vertices of $C_0$, contradicting the choice of $w_0$.

Hence, all edges incident with $w_0$ or $w$ are white and thus belong to the original graph $G$.
We construct an $\{a,b\}$-edge-weighting of $G$ as follows.
Define initial weights $a$ on each edge incident with $w_0$ and extend this weighting to $G$ using~\cref{cor:completion-applied} on $G, H= G-w_0$, and the initial weighting defined on the three white edges incident with $w_0$.
Let $\omega$ denote the resulting $\{a,b\}$-edge-weighting.
For any edge $uv$ in $G-w_0$, the difference $s_\omega(u) - s_\omega(v)$ is not an even integer multiple of $\delta$.

We claim that for any vertex $x \in A(w_0)$, $\da_\omega(x)$ and $\da_\omega(w)$ have different parity.
Since both $x$ and $w$ have $G$-degree $3$, the difference of $b$-degrees is the opposite of the difference of $a$-degrees.
Therefore, we have $s_\omega(x) - s_\omega(w) = (a-b)(\da_\omega(x) - \da_\omega(w))$.
If the $a$-degrees had the same parity, this difference would be an even integer multiple of $\delta$, a contradiction.
This observation implies that the vertices of $A(w_0)$ have $a$-degrees of the same parity.

We now prove that $\omega$ is proper or can be made proper by switching a cycle.
If $\omega$ is proper, we are done; otherwise, the edge $w_0x$ is conflicting for some $x \in A(w_0)$.
Since $s_\omega(w_0) = 3a$ by construction, we have $s_\omega(x) = 3a$.
Then, the parities of the $a$-degrees in $A(w_0)$ are odd, whereas the parity of the $a$-degree of $w$ is even.
Since $d_G(x)=3$ and $s_\omega(x)=3a$, we have $\omega(wx)=a$ and thus $\da_\omega(w)=2$.
Let $y$ be the other vertex of $A(w_0)$ with $\omega(wy) = a$, and let $z$ be the third vertex of $A(w_0)$.
Since $y$ has odd $a$-degree and $\omega(w_0y) = \omega(w y) = a$, we have $s_\omega(y) = 3a$.
Similarly, since $z$ has odd $a$-degree and $\omega(w_0z) = a, \omega(wz) = b$, we have $s_\omega(z) = a + 2b$.

Take an $x-y$ path with internal vertices from $C_0$.
Substitute the corresponding suspended paths of length $3$ for the blue edges in this path to obtain an $x-y$ path in the original graph $G$.
This path, together with the edges $wx$ and $wy$, forms a cycle that avoids $w_0$ and contains two edges incident with each of $x$ and $y$.
After switching this cycle, the edges in $G-w_0$ are still proper by~\cref{cor:completion-applied}.
Furthermore, the edges incident with $w_0$ become proper, since the weighted degree of $w_0$ is still $3a$ and each weighted degree in $A(w_0)$ becomes $a + 2b$.
In conclusion, we obtain a proper $\{a,b\}$-edge-weighting $\omega$ of $G$.

\begin{proof}[Proof of~\cref{thm:main}]
    By~\cref{cor:enfext-case1,prop:enfext-case2,prop:enfext-case3}, either $w_0$ is an enforceable extremum or the exceptional cubic configuration occurs.
    If $w_0$ is an enforceable extremum, then $G$ has a proper $\{a,b\}$-edge-weighting by~\cref{lem:samesign-locminmax,lem:opsign-locmax}.
    In the exceptional configuration, the preceding construction also gives a proper $\{a,b\}$-edge-weighting.
    Both cases contradict the choice of $G$.
\end{proof}

\subsection{Algorithmic aspects}\label{sec:algorithm}

For weights $a,b$ satisfying the hypotheses of~\cref{thm:main}, we first describe the algorithm obtained from the preceding proof when each edge of the input multigraph is given explicitly.
For weights with a rational ratio, we use the normalization in~\cref{sec:proof-completion}; for fixed weights with an irrational ratio, it suffices to construct a proper $\{1,2\}$-edge-weighting by~\cref{lem:irrational-ratio}.
For rational weights given in binary, the normalization takes polynomial time by clearing denominators and dividing by the greatest common divisor.

Among connected bipartite loopless multigraphs with fewer than three vertices, precisely the graphs $K_2^m$ are odd multi-cacti.
For a graph $G$ with a pasted-cycle representation, deleting the endpoints of a bundle $uv$ disconnects $G$ exactly when $uv$ is pasted.
Suppose that $uv$ is pasted and let $g$ be its green-bundle node in the incidence tree $T$.
Each component of $T-g$ gives a subgraph of $G-u-v$ by taking the union of its cycle subgraphs and deleting $u,v$.
These subgraphs cover $G-u-v$ and are pairwise vertex-disjoint.
Each is connected, since its cycle subgraph containing $uv$ remains connected after deletion and the others attach along undeleted green bundles.
Thus, these subgraphs are precisely the components of $G-u-v$.

If $uv$ is not pasted, its unique cycle subgraph remains connected after deletion.
No other cycle subgraph contains both endpoints, since a private green-bundle node has degree one and, for red $uv$, a second cycle node adjacent to both endpoint green-bundle nodes would create a cycle in the incidence tree.
Root the incidence tree at the node of the cycle subgraph containing $uv$.
Adding the other cycle subgraphs in nondecreasing order of the distances of their nodes from the root preserves connectivity, since each loses at most one vertex and attaches along a green bundle with a surviving endpoint.

For graphs with at least three vertices, we test whether $G-u-v$ is disconnected for each edge-bundle $uv$.
If no such bundle is found, any pasted-cycle representation has only one cycle node, so we test whether $G$ is an alternately simple cycle of length $2$ modulo~$4$.
Otherwise, choose such a bundle $uv$ and recursively test the graphs $G[V(Q)\cup\{u,v\}]$ with an extra edge added to $uv$, for each component $Q$ of $G-u-v$.
A representation of $G$ yields a representation of each augmented subgraph by restricting $T$ to the corresponding component of $T-g$ together with $g$, and increasing the multiplicity label of $g$ by one.
Conversely, if all augmented subgraphs admit pasted-cycle representations, then in each representation $uv$ is green by its increased multiplicity and private since deleting $u,v$ leaves $Q$ connected.
Decreasing the corresponding multiplicity labels by one and identifying the private nodes in the disjoint union of the incidence trees gives a representation of $G$.
For $n=|V(G)|\geq3$, the positive integers $|V(H)|-2$ for the subproblems $H$ are additive at each split, so the recursion tree has at most $n-2$ leaves.
Since each internal node has at least two children, there are at most $2(n-2)-1$ nodes in the full recursion tree.
The edge count increases by at most one per level along each branch.

On successful recognition, we return the graph itself if it is $K_2^m$ and a pasted-cycle representation otherwise.
By~\cref{lem:omc-obstruction}, either certificate establishes that no proper weighting exists.

If recognition fails and $G$ has odd order, then~\cref{cor:odd-order} gives a proper weighting.
For even order, we check the properties in~\cref{prop:minimality-prop,prop:redsteps} and apply the corresponding construction if a property fails.
In the reductions of~\cref{prop:minimality-prop}, we test whether the smaller graph $G_0$ is an odd multi-cactus and, if not, recursively construct and extend its proper weighting.
Otherwise, the arguments for~\cref{min:M1,min:M2} exclude the two-vertex case and show that the returned representation of $G_0$, together with the fact that $G$ is not an odd multi-cactus, suffices to apply~\cref{lem:red-modification,lem:green-modification} as appropriate.
We locate the configurations by examining bundles, suspended paths and cycles, and blocks.
The parity completions use the spanning-tree construction in~\cref{lem:parity-subgraph}, and the list selections follow the greedy constructions in~\cref{lem:list-cycle,lem:list-path}.

For the assignment in~\cref{lem:signed-twoholes}, we use dynamic programming, starting with the attainable pair $(0,0)$ and successively adding each admissible pair $(x_i,y_i)$ to the attainable pairs from the preceding step.
At step $i$, we require integers $x_i\in[0,m_i]$ and $y_i\in[0,n_i]$ with $x_i+y_i\equiv\tau_i\pmod2$.
Both coordinates of every attainable pair are bounded by the number of edges, and storing a predecessor for each pair at each step allows us to recover the assignment in polynomial time.
We first seek totals $X,Y$ of the prescribed parities with $X\neq K$ and $Y\neq L$.
If none exist, we recover an assignment with $X=Y=3$ and use the alternative path weighting with both endpoint parity prescriptions reversed, as in the proof of~\cref{lem:opsign-evensuspended}.

The constructions in~\cref{lem:red-modification,lem:green-modification} use the augmented multiplicities from the proof of~\cref{prop:path-replacement} and the list selections above.
Each recursive weighting step uses at most one smaller graph, decreasing the sum of the numbers of vertices and edges.

If no reduction applies, we determine $w_0$ according to the choice rule by testing bundles and computing connected components.
Since $A$ is independent in $J$, the sets $E_J(x)$ for $x\in A$ are disjoint.
The spanning tree of~\cref{cor:frank-leaving} can therefore be found by polynomial-time matroid intersection~\cite{Frank2011} on the graphic matroid of $J$ and the partition matroid with classes $E_J(x)$ for $x\in A$ and the remaining edges, imposing bounds $d_J(x)-h(x)-1$ on the former classes and no restriction on the latter.
If any bound is negative, no such tree exists; otherwise we compute a maximum-cardinality common independent set, which gives such a tree exactly when it has size $|V(J)|-1$.
If we obtain such a spanning tree $T$, the edges of $E(J)\setminus E(T)$ with an endpoint in $A$ form a set $F_R$ satisfying $d_{F_R}(x)=d_J(x)-d_T(x)\geq h(x)+1$ for each $x\in A$, while $J-F_R$ contains $T$.
The conversion described above then gives $F$.
The proofs of~\cref{cor:enfext-case1,prop:enfext-case2,prop:enfext-case3} show that the absence of such a spanning tree can occur only in the exceptional cubic configuration, for which the preceding cycle-switching construction gives a proper weighting.
After normalization or the replacement by $\{1,2\}$, the number of operations in these constructions is bounded polynomially in the numbers of vertices and edges, independently of the numerical values of the integer weights.
Weighted degrees, fixed contributions, and list entries are sums of the integer weights with at most as many terms as there are edges in the graph under consideration.
In the proof of~\cref{lem:completion-general}, we may choose $\alpha_Q=\theta(v)$ for any $v\in V(Q)$, so the integrality and parity tests use only integer arithmetic on numerators and the common denominator $b-a$.
All intermediate numbers have polynomial bit length in the encoding lengths of the graph and the integer weights.
Thus each step takes polynomial time in these encoding lengths, and there are only polynomially many recursive steps.

The algorithm also works when the multiplicities of the edge-bundles are given in binary and the output weighting is specified by the number of edges receiving weight $b$ in each bundle, given in binary.
Let $n=|V(G)|$ and truncate each multiplicity $m$ to $\min\{m,n\}$, obtaining a graph $G'$ with $O(n^3)$ edges.
By~\cref{prop:models}, truncation preserves membership in the odd-multi-cactus family, since it preserves the underlying simple graph and the distinction between simple and multiple bundles.
If recognition succeeds on $G'$, restore the original multiplicities in the returned certificate in binary.
Otherwise, let $\alpha,\beta$ denote the integer weights used by the explicit-edge algorithm.
Construct a proper $\{\alpha,\beta\}$-edge-weighting of $G'$ and restore the original multiplicities one bundle at a time.
When restoring a bundle $uv$ of multiplicity $m>n$, consider the $n$ distinct totals $m\alpha+k(\beta-\alpha)$ for $k\in[0,n-1]$.
Each choice changes both endpoint sums by the same amount, so $uv$ remains proper.
Since $G$ is bipartite, at most $n-2$ other bundles are incident with $u$ or $v$, and each forbids at most one total.
Hence a choice preserving properness exists at every step.
For each of the $O(n^2)$ bundles, we test at most $n$ totals using integers of polynomial bit length.
Thus the running time is polynomial in the binary input length for fixed weights or rational weights given as input.
For a bundle of multiplicity $m$ with $k$ edges receiving $\beta$, we return $k$ if $\beta$ corresponds to $b$ and $m-k$ otherwise.
In the irrational-ratio case, $1$ corresponds to $a$ and $2$ to $b$.

\section{Universal 3-edge-weightability}\label{sec:universal}

We now prove~\cref{thm:universal}, establishing universal $3$-edge-weightability for connected bipartite multigraphs of order at least three.
By~\cref{thm:main}, only odd multi-cacti and triples $\{-c,0,c\}$ with $c>0$ require separate arguments.
We first prove a stronger list-weighting result for odd multi-cacti.

\begin{lemma}\label{lem:omc-list}
Let $G$ be an odd multi-cactus with at least three vertices and let $L(e)\subseteq\R$ be a set of three weights for every $e\in E(G)$.
Then $G$ admits a proper edge-weighting $\omega$ with $\omega(e)\in L(e)$ for every $e\in E(G)$.
\end{lemma}

\begin{proof}
Fix a pasted-cycle representation of $G$ by~\cref{prop:models}.
We first weight the red edges to make the green bundles proper, and then weight the green bundles to make the red edges proper.
As in the proof of~\cref{lem:cycle-bundle}, the total of a green bundle $uv$ cancels in $s(u)-s(v)$, leaving only red weights.
Each red edge occurs in exactly two such differences with coefficient $1$ or $-1$, since its endpoints belong to distinct green bundles.
Weight the red edges in any order, requiring each difference to be non-zero when its last red edge is weighted.
At most two weights are forbidden at each step, so a choice from the three-element list exists.

Obtain a graph $H$ from $G$ by identifying the two endpoints of each green edge and deleting all green edges.
The edges of $H$ correspond to the red edges of $G$.
By the attachment order in the proof of~\cref{prop:models}, $H$ is obtained by successively attaching cycles of length at least three along single vertices.
Thus $H$ is simple and $2$-degenerate, as seen by deleting new vertices in reverse attachment order.
Each green bundle admits at least three totals by fixing an admissible weight for all but one of its edges and varying the weight of the remaining edge over its own list.
Process the vertices of $H$ in an order with at most two earlier neighbors per vertex, choosing a total for each green bundle to make the red edges to earlier vertices proper.
At most two totals are forbidden at each step, and these choices do not affect properness across green bundles.
\end{proof}

\begin{proof}[Proof of~\cref{thm:universal}]
Let $G$ be as in~\cref{thm:universal} and let $W\subseteq\R$ have three elements.
By~\cref{cor:odd-order}, we may assume that $G$ has even order.
If $G$ is an odd multi-cactus, apply~\cref{lem:omc-list} with $L(e)=W$ for all edges $e$.
Otherwise, \cref{thm:main} applies unless $W=\{-c,0,c\}$ for some $c>0$, since every other triple contains two non-zero weights with non-zero sum.
By scaling, assume $W=\{-1,0,1\}$.

Let $A,B$ be the bipartition classes of $G$ and fix a spanning tree $T$ of $G$.
If $|A|$ is even, use~\cref{lem:parity-subgraph} in $T$ to find an edge set whose odd-degree vertices are precisely $A$.
Weights $1$ on this set and $0$ elsewhere give a proper weighting by parity.
Suppose both classes have odd cardinality, with $|A|\geq3$ without loss of generality.
Then one exception to the parity prescription is necessary; we make that vertex a strict weighted local minimum.
Fix $v\in A$ and choose $F\subseteq E(T)$ with odd-degree vertices $A\setminus\{v\}$ by~\cref{lem:parity-subgraph}.
In particular, $F\neq\emptyset$.
If $d_F(v)>0$, put $w=v$; otherwise, choose a vertex $w$ of positive $F$-degree closest to $v$ in $T$ and add the edges of the $v$--$w$ path in $T$ to $F$.
The added path, if any, has no edge in the original $F$ and changes parities only at $v$ and $w$.
We have $d_F(w)\geq2$, since its $F$-degree is positive and even if $w=v$, and otherwise the added path increases its positive $F$-degree by one.
For every vertex $u\neq w$, the degree $d_F(u)$ is odd if and only if $u\in A$.

Assign weight $-1$ to the edges of $F$ incident with $w$, $1$ to the other edges of $F$, and $0$ elsewhere.
Since $s(u)\equiv d_F(u)\pmod2$ for all vertices $u$, all edges not incident with $w$ are proper.
Each neighbor $u$ of $w$ has $s(u)\geq-1$, since $F\subseteq E(T)$ contains at most one edge between $w$ and $u$.
As $s(w)=-d_F(w)\leq-2$, the edges incident with $w$ are proper as well.
\end{proof}

For rational weights in binary, the proof and the algorithms of~\cref{sec:algorithm} yield a polynomial-time algorithm.
For binary multiplicities, use the same truncation and restoration with any two weights, avoiding at most $|V(G)|-2$ totals per bundle.
The output gives binary counts per weight and bundle.

\section*{Acknowledgments}

We thank Carsten Thomassen for helpful correspondence concerning the proof of Theorem~5.1 in~\cite{thomassenwuzhang2016flows}.

\smallskip
\noindent This research has been implemented with the support provided by the Ministry of Innovation and Technology of Hungary from the National Research, Development and Innovation Fund, financed under the ELTE TKP 2021-NKTA-62 funding scheme.

\section*{Declarations}

Writefull, Gemini, and ChatGPT were used for language editing.
The authors take full responsibility for the content of the manuscript.

\smallskip
\noindent The authors declare that they have no conflicts of interest.

\AddToHook{cmd/thebibliography/after}{\setlength{\itemsep}{.1pt}}
\bibliographystyle{plain}
\bibliography{bibliography}

\begin{thebibliography}{10}

\bibitem{bartnicki2009weight}
Tomasz Bartnicki, Jaros{\l}aw Grytczuk, and Stanis{\l}aw Niwczyk.
\newblock Weight choosability of graphs.
\newblock {\em Journal of Graph Theory}, 60(3):242--256, 2009.

\bibitem{bensmail2022oddab}
Julien Bensmail, Fionn {Mc Inerney}, and Kasper~Szabo Lyngsie.
\newblock On {$\{a,b\}$}-edge-weightings of bipartite graphs with odd {$a,b$}.
\newblock {\em Discussiones Mathematicae Graph Theory}, 42(1):159--185, 2022.

\bibitem{chang2011graphs}
Gerard~J. Chang, Changhong Lu, Jiaojiao Wu, and Qinglin Yu.
\newblock Vertex-coloring edge-weightings of graphs.
\newblock {\em Taiwanese Journal of Mathematics}, 15(4):1807--1813, 2011.

\bibitem{deng2026regular}
Kecai Deng.
\newblock Regular graphs are universally {$3$}-edge-weightable, 2026.
\newblock arXiv:2602.06659v2.

\bibitem{dengqiu2027uniformspan}
Kecai Deng and Hongyuan Qiu.
\newblock Every graph is uniform-span {$(2,2)$}-choosable: Beyond the {1-2}
  conjecture.
\newblock {\em Journal of Combinatorial Theory, Series B}, 182:40--49, 2027.
\newblock arXiv:2506.14253v8.

\bibitem{dudekwajc2011complexity}
Andrzej Dudek and David Wajc.
\newblock On the complexity of vertex-coloring edge-weightings.
\newblock {\em Discrete Mathematics \& Theoretical Computer Science},
  13(3):45--50, 2011.

\bibitem{Frank2011}
Andr{\'a}s Frank.
\newblock {\em Connections in Combinatorial Optimization}, volume~38 of {\em
  Oxford Lecture Series in Mathematics and Its Applications}.
\newblock Oxford University Press, Oxford, 2011.

\bibitem{karonski2004edge}
Micha{\l} Karo{\'n}ski, Tomasz {\L}uczak, and Andrew Thomason.
\newblock Edge weights and vertex colours.
\newblock {\em Journal of Combinatorial Theory, Series B}, 91(1):151--157,
  2004.

\bibitem{keusch2024solution123}
Ralph Keusch.
\newblock A solution to the {1-2-3} conjecture.
\newblock {\em Journal of Combinatorial Theory, Series B}, 166:183--202, 2024.

\bibitem{khatirinejad2012twoedgeweights}
Mahdad Khatirinejad, Reza Naserasr, Mike Newman, Ben Seamone, and Brett
  Stevens.
\newblock Vertex-colouring edge-weightings with two edge weights.
\newblock {\em Discrete Mathematics \& Theoretical Computer Science},
  14(1):1--20, 2012.

\bibitem{lu2016bipartite}
Hongliang Lu.
\newblock Vertex-coloring edge-weighting of bipartite graphs with two edge
  weights.
\newblock {\em Discrete Mathematics \& Theoretical Computer Science},
  17(3):1--12, 2016.

\bibitem{luyuzhang2011twoedgeweighting}
Hongliang Lu, Qinglin Yu, and Cun-Quan Zhang.
\newblock Vertex-coloring {$2$}-edge-weighting of graphs.
\newblock {\em European Journal of Combinatorics}, 32(1):21--27, 2011.

\bibitem{lyngsie2018nsd01}
Kasper~Szabo Lyngsie.
\newblock On neighbour sum-distinguishing {$\{0,1\}$}-weightings of bipartite
  graphs.
\newblock {\em Discrete Mathematics \& Theoretical Computer Science}, 20(1),
  2018.
\newblock Article 21.

\bibitem{madarasi2026trees}
P{\'e}ter Madarasi.
\newblock Proper $\{a,b\}$-edge-weightings of trees.
\newblock {\em Discrete Mathematics}, 350(4):115539, 2027.

\bibitem{madarasi2026complexity}
P{\'e}ter Madarasi and M{\'a}t{\'e} Simon.
\newblock Complexity, approximation, and extension of proper
  $\{a,b\}$-edge-weightings, 2026.
\newblock arXiv:2609.31205.

\bibitem{thomassenwuzhang2016flows}
Carsten Thomassen, Yezhou Wu, and Cun-Quan Zhang.
\newblock The {$3$}-flow conjecture, factors modulo {$k$}, and the
  {1-2-3}-conjecture.
\newblock {\em Journal of Combinatorial Theory, Series B}, 121:308--325, 2016.

\end{thebibliography}

\end{document}